\RequirePackage{fix-cm}
\documentclass{article}
\usepackage{graphicx,amsfonts,amsmath}

\usepackage{amssymb}
\usepackage{hyperref}

\usepackage{cite}

\usepackage{color}

\usepackage{url}

\usepackage{subfigure}
\usepackage{enumitem}

\usepackage[a4paper, total={6in, 8in}]{geometry}

\newcommand{\inner}[2]{\langle#1,#2\rangle}

\renewcommand{\Re}{\mathbb{R}}

\newcommand{\E}{\mathbb{E}}
\newcommand{\Q}{\mathcal{Q}}

\newcommand{\exclude}[1]{}

\newcommand{\col}[1]{\left\{#1\right\}}

\newcommand{\scr}[1]{\mathcal{#1}} 

\newcommand{\fc}[2]{: #1 \rightarrow #2}

\makeatletter
\newcommand\tabsize{\@setfontsize\tinyv{10pt}{12}}
\makeatother

\makeatletter
\newcommand\tabsizesm{\@setfontsize\tinyv{8pt}{10}}
\makeatother

\newtheorem{remark}{Remark}[section]

\begin{document}

\title{On conditional cuts for Stochastic Dual Dynamic
Programming}

\title{On conditional cuts for Stochastic Dual Dynamic Programming}        

\author{W. van Ackooij    \footnote{ EDF R\&D, \sf \href{mailto:wim.van-ackooij at edf.fr}{wim.van-ackooij at edf.fr} }   \and
        X.~Warin  \footnote{EDF R\&D \& FiME \sf \href{mailto:xavier.warin at edf.fr}{xavier.warin at edf.fr}}
}


\maketitle

\begin{abstract}
Multi stage stochastic programs arise in many applications from
engineering whenever a set of inventories or stocks has to be
valued. Such is the case in seasonal storage valuation of a set of
cascaded reservoir chains in hydro management. A popular method is
Stochastic Dual Dynamic Programming (SDDP), especially when the
dimensionality of the problem is large and Dynamic programming is
no longer an option.  The usual assumption of SDDP is that
uncertainty is stage-wise independent, which is highly restrictive
from a practical viewpoint. When possible, the usual remedy is to
increase the state-space to account for some degree of dependency.
In applications this may not be possible or it may increase the
state space by too much. In this paper we present an alternative
based on keeping a functional dependency in the SDDP - cuts
related to the conditional expectations in the dynamic programming
equations. Our method is based on popular methodology in
mathematical finance, where it has progressively replaced scenario
trees due to superior numerical performance. We demonstrate the
interest of combining this way of handling dependency in
uncertainty and SDDP on a set of numerical examples. Our method is
readily available in the open source software package StOpt.


\end{abstract}

\section{Introduction}\label{sec:introduction}

Dealing with uncertainty is vital in many real-life applications.
An interesting way to formalize such a setting is a multistage
stochastic program wherein one also accounts for the possibility
of acting on past observed uncertainty. These models are quite
popular in areas such energy
\cite{Goel_Grossmann_2004,Pereira_Granville_Fampa_Dix_Barroso_2005,Shapiro_Tekaya_daCosta_Soares_2013,Rebennack_2016,deMatos_Morton_Finardi_2016},
transportation
\cite{Fhoula_Hajji_Rekik_2013,Herer_Tzur_Yucesan_2006} and finance
\cite{Pflug_Romisch_2007,Dupacova_2009,Dupacova_Polivka_2009}. The
usual underlying assumption is that uncertainty can somehow be
presented or approximated using a scenario tree. The resulting
mathematical programming problem is generally large-scale and
non-trivial or impossible to solve using a monolithic method. The
use of specialized algorithms that employ decomposition techniques
(and very often sampling) appear crucial for an efficient
numerical resolution. In this family two popular approaches are
\emph{Nested Decomposition} -- ND -- of \cite{Birge_1985b} and
\emph{Stochastic Dual Dynamic Programming} -- SDDP -- of
\cite{Pereira_Pinto_1991}. Mathematically speaking, both methods
are close cousins, in that they approximate the cost-to-go
functions (resulting from dynamic programming) using piecewise
linear approximations, defined by cutting planes computed by
solving linear programs in a so-called \emph{backward step}.

The difference between the ND and SDDP algorithm resides in how
uncertainty is handled when trying to establish an upper bound on
the optimal value. The ND algorithm will use the entire scenario
tree, whereas the SDDP algorithm will employ sampling procedures
in order to achieve this. The ND algorithm therefore does not
require any particular assumption on the scenario tree, but as a
consequence is typically only applied to multistage stochastic
problems of moderate size (with some hundred or a few thousand
scenarios). In order to mitigate the effect that only ``smaller''
trees can be handled, methods for constructing representative but
small scenario trees have therefore received significant
attention. Let us mention the pioneering work
\cite{Dupacova_Growe-Kuska_Romisch_2003} on two-stage programs and
some subsequent extensions to the multistage case: e.g.,
\cite{Oliveira_Sagastizabal_Jardim-Penna_Maceira_Damazio_2010,Pflug_2010,Heitsch_Romisch_2011,Pflug_Pichler_2012,Kovacevic_Pichler_2015}.
Other ideas to reduce the computational burden rely on combining
adaptive partitioning of scenarios with regularization techniques
from convex optimization, e.g.,
\cite{Song_Luedtke_2015,vanAckooij_Oliveira_Song_2018}. These
ideas can also be extended to the multistage setting as done
recently in \cite{vanAckooij_Oliveira_Song_2019}.

It is therefore generally acknowledged that SDDP can handle larger
scenario trees by combining decomposition and sampling,
\cite{Pereira_Pinto_1991}. This method is highly popular as
illustrated by the numerous applied papers, e.g.,
\cite{Oliveira_Sagastizabal_Jardim-Penna_Maceira_Damazio_2010,Shapiro_Tekaya_daCosta_Soares_2013,Herer_Tzur_Yucesan_2006,Philpott_deMatos_2012,Goel_Grossmann_2004,Rodriguez_Oliveira_finardi_2017}.
The method however requires the assumption that the underlying
stochastic process is stagewise independent (see also
\cite{Chen_Powell_1999,Hindsberger_Philpott_2001,Donohue_Birge_2006}
for other methods employing sampling). Combined with a possibility
to share cuts among different nodes of the scenario tree, the
optimization strategies proposed in these references mitigate the
curse of dimensionality further. However, the assumption that
uncertainty is stagewise independent is quite a strong assumption
when familiar with the actual observation of data. Whenever the
underlying stochastic process is Markovian, it is generally
suggested to increase the dimension of the state vector in order
to revert back to a stagewise independence assumption (c.f., the
discussion in \cite{Shapiro_2011}). As observed in
\cite{deQueiroz_Morton_2013} such an increase may be detrimental
to computational efficiency especially if the increase in the
dimension of the state vector is significant.
%
%
The authors \cite{deQueiroz_Morton_2013} suggest a way to
partially mitigate this effect. In this paper we suggest another
approach for efficiently estimating the conditional expectations
appearing in the dynamic programming equations. The suggested
method does not require increasing the size of the state vector
whenever the underlying stochastic process is not stagewise
independent. Our approach of estimating the conditional
expectation is based on the observation that such expectations are
orthogonal projections onto an appropriate functional space. The
method suggested here is very popular in mathematical finance
since the first work of Tsilikis and Van Roy
\cite{Tsilikis_VanRoy_1999} and the closely related Longstaff
Schwarz method \cite{Longstaff_Schwartz_2001}. The latter has
become the most used method to deal with optimal stopping problem,
i.e., the optimization problem related to computing the optimal
stopping time. The optimal stopping time is the best moment to
exercise a given option. In the banking system, this method has
become the reference method due to its easy implementation, its
efficiency in moderate dimensions, and the easiness to calculate
 the sensibility of the optimal value (the ``greeks'').
In Energy market, quants use this method to valuate and hedge
their portfolio (see \cite{Warin_2012} for an example of a gas
storage hedging simulation).

In a recent paper, Bouchard and Warin \cite{Bouchard_Warin_2012}
developed a variant of this method based on linear regression on
adaptive meshes: they showed that this method was superior to the
original Longstaff-Schwarz method based on global polynomial
regression by avoiding oscillations of the regressor (the Runge
effect). Indeed, the method was compared to different methods to
price American options and to estimate the sensitivity of the
options to the initial prices (the delta of the option) in
dimension 1 to 6. In dimension 3 or above, the proposed regression
method appeared to be clearly superior to optimal quantization
\cite{Pages_Printems_2003,Pages_Pham_Printemps_2004,Pages_Printemps_2008,Pages_Printemps_2005}
and the Malliavin approach
\cite{Fournie_Lasry_Lebuchoux_Lions_Touzi_1999,Fournie_Lasry_Lebuchoux_Lions_2001}.
Finally let us mention that the here suggested regression methods
are also used to approximate conditional expectations in numerical
schemes \cite{Bouchard_Touzi_2004} used to solve Backward
Stochastic Differential Equations (BSDE). The latter BSDEs provide
a way to solve quasi-linear equations: the first papers
\cite{Gobet_Lemor_Warin_2005,Lemor_Gobet_Warin_2006} have given
birth to many follow up papers on the topic (see for example the
references in the recent work of \cite{Gobet_Plamen_2016}). One
can also use these building block to design schemes to solve full
nonlinear equations \cite{Fahim_Touzi_Warin_2011}.

To conclude the introduction, let us briefly mention that the
mathematical properties of the SDDP algorithm have been
extensively investigated. A first formal proof of almost-sure
convergence of multistage sampling algorithms akin to SDDP is due
to Chen and Powell \cite{Chen_Powell_1999}. This proof was
extended by Linowsky and Philpott to cover the SDDP algorithm and
other sampling-based methods in~\cite{Linowsky_Philpott_2005}. In
fact both papers rely on a tacit assumption, formally worked out
in \cite{Philpott_Guan_2008}. Convergence in the general convex
risk-neutral case is worked out in
\cite{Girardeau_Leclere_Philpott_2015} and also studied in the
risk-averse situation in \cite{Guigues_2016}. Recent publications
have addressed the issue of incorporating risk-averse measures
into the SDDP algorithm:
\cite{Guigues_Romisch_2012,Shapiro_2011,Philpott_deMatos_2012,Guigues_Sagastizabal_2012,Shapiro_Tekaya_daCosta_Soares_2013,
Dupacova_Valclav_2015,Homem-de-Mello_Pagnoncelli_2016}. Polyhedral
risk measures were studied in \cite{Eichhorn_Romisch_2005} and
\cite{Guigues_Romisch_2012}. As mentioned
in~\cite{Shapiro_Tekaya_daCosta_Soares_2013}, theoretical
foundations for a risk averse approach based on conditional risk
mappings were developed in \cite{Ruszczynski_Shapiro_2006b}.
Extended polyhedral risk measures, also allowing for dynamic
programming equations, can also be employed and where developed in
\cite{Guigues_Romisch_2012}. It is shown is \cite{Shapiro_2011}
how to incorporate convex combinations of the expectation and
\emph{Average Value-at-Risk} into the SDDP algorithm. Numerical
experiments on this idea have been reported in many publications;
see for instance \cite{Philpott_deMatos_2012} and
\cite{Shapiro_Tekaya_daCosta_Soares_2013}.

This paper is organized as follows. Section \ref{sec:decomp}
presents the general structure of multi-stage stochastic linear
programs and show how the key difficulty resides in computing
conditional expectations. We present the general idea of
approximating conditional expectations in Section
\ref{sec:condcut} and our algorithm in section \ref{algo:genSDDP}.
Finally Section \ref{sec:numerical} contains a series of numerical
experiments showing numerical consistency of the suggested method.
We also compare the method with other variants, wherein it appears
that the suggested approach is preferable. The paper ends with
concluding remarks and research perspectives.

\section{Preliminaries on multistage stochastic linear programs and decomposition}
\label{sec:decomp}

Multistage stochastic programs explicitly model a series of
decisions interspaced with the partial observation of uncertainty.
If the given set of possible realizations of the underlying
stochastic process is discrete, uncertainty can be represented by
a scenario tree. Multistage stochastic linear programs (MSLPs) are
a special case of this class wherein the underlying modelling
structure is linear/affine. Hence, on a scenario tree the
resulting model is a very large linear program. In this section we
will provide an overview of two popular decomposition approaches
for MSLPs, namely nested decomposition and stochastic dual dynamic
programming. From an abstract viewpoint coming from non-linear
non-smooth optimization, both methods are variants of Kelley's
cutting plane method \cite{Kelley_1960}.

Consider the following multistage stochastic linear program:
\begin{equation}\label{mslp}
\min_{\overset{x_1\in X_1}{A_1x_1=b_1}} c_1^\top  x_1+
\E_{|\xi_1}\left[ \min_{\overset{x_2\in X_2}{B_2x_1 +
A_2x_2=b_2}}c_2^\top x_2 + \E_{|\xi_{[2]}}
\left[\cdots+\E_{|\xi_{[T-1]}}   [ \min_{\overset{x_T\in
X_T}{B_Tx_{T-1} + A_Tx_T=b_T}}c_T^\top  x_T ]\right]\right] \,,
\end{equation}
where some of (or all) data $\xi=(c_t,B_t,A_t,b_t)$ can be subject
to uncertainty for $t=2,\ldots,T$. The conditional expectation
$\E_{|\xi_{[t]}}[\,\cdot\,]$ is taken with respect to the
filtration generated by the history, up to time $t$, of the random
vector $\xi_t\in \Xi_t \subseteq \Re^{m_t}$, defining the
stochastic process $\{\xi_{t}\}_{t=1}^T$. This history, up to time
$t$, will be denoted with $\xi_{[t]} =(\xi_1,\ldots,\xi_t)$.
Furthermore, $X_t\neq \emptyset$, $t=1,\ldots,T$, are polyhedral
convex sets that do not depend on the random parameters, which we
denote by $X_t := \{x_t\in \Re_+^{n_t} \mid D_tx_t = d_t\}$, for
an appropriate matrix $D_t$ and vector $d_t$.

For numerical tractability, it is typically assumed that the number $N$ of
realizations (scenarios) of the data process is finite, i.e.,
support sets $\Xi_t$ ($t=1,\ldots,T$) have finite cardinality.
This is, for instance, the case in which~\eqref{mslp} is a SAA
approximation of a more general MSLP problem (having continuous
probability distribution). For a discussion of the relation of
such a problem with the underlying true problem relying on a
continuous distribution we refer to \cite{Shapiro_2011}. In our numerical experiments we will also make this assumption.

Under these assumptions, the dynamic programming equations for
problem \eqref{mslp} take the form
\begin{equation}\label{cost2go}
Q_t(x_{t-1},\xi_{[t-1]},\xi_{t}):= \left\{
\begin{array}{llll}
\displaystyle\min_{x_t\in\Re^{n_t}} & c_t^\top  x_t + \Q_{t+1}(x_t)(\xi_{[t]})\\
\mbox{s.t.} &A_tx_t=b_t-B_tx_{t-1} \\
 &   x_t \in X_t\,,
\end{array}
\right.
\end{equation}
where
\begin{equation}\label{recurse}
\Q_{t+1}(x_t)(\xi_{[t]}) :=
\E_{|\xi_{[t]}}[Q_{t+1}(x_t,\xi_{[t]},\xi_{t+1})]\,,\quad
\mbox{for }\; t=T-1,\ldots 1\,,
\end{equation}
and $\Q_{T+1}(x)\equiv 0$ by definition. The first-stage problem
becomes
\begin{equation}\label{1stagepbm}
\left\{
\begin{array}{llll}
\displaystyle \min_{x_1\in\Re^{n_1}}& c_1^\top x_1 + \Q_2(x_1)(\xi_{[1]})\\
\mbox{s.t.}& A_1x_1=b_1\\
           & x_1 \in X_1\,.
           \end{array}
           \right.
\end{equation}
An implementable policy for \eqref{mslp} is a collection of
functions $\bar x_t=\bar x_t(\xi_{[t]})$, $t=1,\ldots,T$. Such a
policy gives a decision rule at every stage $t$ of the problem
based on a realization of the data process up to time $t$. A
policy is feasible for problem~\eqref{mslp} if it satisfies all
the constraints for every stage $t$.

As in \cite{Shapiro_2011}, we  assume that the \emph{cost-to-go}
functions $\Q_t$ are finite valued, in particular we assume
\emph{relatively complete recourse}. Since the number of scenarios
is finite, the cost-to-go functions are convex piecewise linear
functions \cite[Chap. 3]{Shapiro_Dentcheva_Ruszczynski_2009}.

We note that the stagewise independence assumption (e.g., as in
\cite{Pereira_Pinto_1991}) simplifies \eqref{recurse} to
$\Q_{t+1}(x_t)(\xi_{[t]}) = \Q_{t+1}(x_t) =
\E[Q_{t+1}(x_t,\xi_{[t]},\xi_{t+1})]$ since the random vector
$\xi_{t+1}$ is independent of its history $\xi_{[t]}
=(\xi_1,\ldots,\xi_t)$. As pointed out in \cite{Shapiro_2011}, in
some cases stagewise dependence can be dealt with by adding state
variables to the model. We will provide an alternative to that
general methodology. We also care to note, and this will be a
subsequent assumption, that whenever the data process is
Markovian, than then $\Q_{t+1}(x_t)(\xi_{[t]})$ depends on $\xi_t$
alone, i.e.,
$\Q_{t+1}(x_t)(\xi_{[t]})\equiv\Q_{t+1}(x_t)(\xi_{t})$.

First, let us briefly mention the main steps of ND and SDDP.

\subsection{Decomposition}

As already mentioned, two very important decomposition techniques
for solving multistage stochastic linear programs are Nested
Decomposition (see \cite{Ruszczynski_Swietanowski_1997}), and the
SDDP algorithm (see \cite{Pereira_Pinto_1991}). Both methods have
two main steps:
\begin{itemize}
\item \emph{Forward step}, that goes from stage $t = 1$ up to $t =
T$ solving subproblems to define feasible policies $\bar
x_t(\xi_{[t]})$. In this step an (estimated) upper bound
$\overline{z}$  for the optimal value is determined.

\vspace{.3cm} \item \emph{Backward step}, that comes from stage $t
= T$ up to $t = 1$ solving subproblems to compute linearizations
that improve the cutting-plane approximations for the cost-to-go
functions $\Q_t$. In this step a lower bound $\underline z$ is
obtained.
\end{itemize}
Below we discuss these steps, starting with the backward one.

\paragraph{Backward step.} Let $\bar x_t=\bar x_t(\xi_{[t]})$
be a trial decision at stage $t = 1,\ldots, T-1$, and
$\check{\Q}_t$ be a current approximation of the cost-to-go
function $\Q_t$, $t = 2,\ldots, T$, given by the maximum of a
collection of cutting planes. At stage $t = T$ the following
problem is solved
\begin{equation}\label{cpT}
\underline{Q}_T(\bar x_{T-1},\xi_{[T-1]},\xi_T)= \left\{
\begin{array}{llll}
\displaystyle\min_{x_T\in\Re^{n_T}} & c_T^\top  x_T\\
\mbox{s.t.} & A_Tx_T=b_T - B_T\bar x_{T-1} \\
& x_T\in X_T
\end{array}
\right.
\end{equation}
for all $\xi_T=(c_T,B_T,A_T,b_T) \in \Xi_T$. Let $\bar \pi_T= \bar
\pi_T(\xi_T)$ be an optimal dual solution of problem \eqref{cpT}.
Then $\alpha_T(\xi_{T-1}) := \E_{|\xi_{T-1}}[b_T^\top \bar \pi_T]$
and $ \beta_T(\xi_{T-1}):= -\E_{|\xi_{T-1}}[B_T^\top \bar
\pi_T]\in \partial \Q_T(\bar x_{T-1})(\xi_{T-1})$ define the
linearization
\begin{align*}
q_T(x_{T-1})(\xi_{T-1}) &:= \beta_T^{\top}(\xi_{T-1}) x_{T-1} +
\alpha_T(\xi_{T-1}) \\
&= \quad \Q_T(\bar x_{T-1})(\xi_{T-1}) +
\inner{\beta_T(\xi_{T-1})}{x_{T-1} - \bar x_{T-1}}\,,
\end{align*}
satisfying
\[ 
\Q_{T}(x_{T-1}) \geq q_T (x_{T-1}) \quad \forall \; x_{T-1}\,,
\]
 and $\Q_{T}(\bar x_{T-1})(\xi_{T-1}) = q_T(\bar x_{T-1})(\xi_{T-1})$, i.e., $q_T$ is a supporting
plane for $\Q_{T}$ (both functions of $\xi_{T-1}$). This
linearization is added to the collection of supporting planes of
$\Q_T$: $\check \Q_T$ is replaced by
$\check\Q_T(x_{T-1})(\xi_{T-1}):= \max\{\check
\Q_T(x_{T-1})(\xi_{T-1}),q_T(x_{T-1})(\xi_{T-1})\} $. In other
words, the cutting-plane approximation $\check \Q_T$ is
constructed from a collection $J_T$ of linearizations:
\[
\check\Q_T(x_{T-1})(\xi_{T-1}) = \max_{j \in
J_T}\{\beta_{T}^{j\top}(\xi_{T-1}) x_{T-1} +
\alpha_{T}^j(\xi_{T-1}) \}\,.
\]

By starting with a model approximating the cost-to-go function
from below (e.g. $\Q_t\equiv -\infty$ for all stages), the
cutting-plane updating strategy will ensure that $\check \Q_T\leq
\Q_T$. The updated model $\check \Q_T$ is then used at stage
$T-1$, and the following problem needs to be solved for all
$t=T-1,\ldots, 2$:
\begin{align}
\underline{Q}_{t}(\bar x_{t-1},\xi_{[t-1]},\xi_{t}) &=\left\{
\begin{array}{llll}
\displaystyle \min_{x_{t}\in\Re^{n_t}}& c_{t}^\top  x_{t} + \check \Q_{t+1}(x_{t})(\xi_{t})\\
\mbox{s.t.}& A_{t}x_{t}=b_{t} - B_{t}\bar x_{t-1}\\
& x_{t}\in X_{t}
\end{array}
\right. \nonumber \\
& \equiv \left\{
\begin{array}{llll}
\displaystyle \min_{x_{t},r_{t+1}\in\Re^{n_t}\times\Re}& c_{t}^\top  x_{t} + r_{t+1}\\
\mbox{s.t.}&  A_{t}x_{t}=b_{t} - B_{t}\bar x_{t-1}\\
 & \beta_{t+1}^{j\top}(\xi_{t}) x_{t} + \alpha_{t+1}^j(\xi_{t}) \leq r_{t+1},\quad j \in J_{t+1}\\
 &  x_{t}\in X_{t}\,.
\end{array}
\right. \label{cpt-1}
\end{align}
Let $\bar \pi_t=\bar \pi_t(\xi_t)$ and $\bar \rho_t^{j} = \bar
\rho_t^j(\xi_t)$ be optimal Lagrange multipliers associated with
the constraints $A_{t}x_{t}=b_{t} - B_{t}\bar x_{t-1}$ and
$\beta_{t+1}^{j\top}(\xi_{t-1}) x_{t} + \alpha_{t+1}^j(\xi_{t-1})
\leq r_{t+1}$, respectively. Then the linearization
\[
q_t(x_{t-1})(\xi_{t-1}):= \beta_t^{\top}(\xi_{t-1}) x_{t-1} +
\alpha_t(\xi_{t-1}) = \E_{|\xi_{t-1}}[\underline{Q}_{t}(\bar
x_{t-1},\xi_{[t-1]},\xi_{t})] + \inner{\beta_t(\xi_{t-1})}{x_{t-1}
- \bar x_{t-1}}
\]
of $\Q_t$ is constructed with
 \begin{equation}\label{cutcoeff}
 \alpha_t(\xi_{t-1}):= \E_{|\xi_{t-1}}[b_t^\top \bar \pi_t + \sum_{j\in J_{t+1}} \alpha_{t+1}^j(\xi_{t}) \bar \rho_t^j] \quad\mbox{and}\quad
 \beta_t(\xi_{t-1}):= -\E_{|\xi_{t-1}}[B_t^\top \bar \pi_t]
 \end{equation} and satisfies
$ \Q_{t}(x_{t-1})(\xi_{t-1}) \geq q_t (x_{t-1})(\xi_{t-1})\;\;
\forall \; x_{t-1}$.

Once the above linearization is computed, the cutting-plane model
at stage $t$ is updated: $\check\Q_t(x_{t-1})(\xi_{t-1})
=\max\{\check\Q_t(x_{t-1})(\xi_{t-1}),q_t(x_{t-1})(\xi_{t-1})\}$\,.
Since $\check \Q_{t+1}$ can be a rough approximation (at early
iterations) of $\Q_{t+1}$, the linearization $q_t$ is not
necessarily a supporting plane (but a cutting plane) for
$\Q_{t+1}$: the inequality $q_t\leq \Q_t$ might be strict for all
$x_{t-1}$ feasible (at the first iterations).

At the first stage, the following LP is solved
\begin{align}
\underline{z} &=\left\{
\begin{array}{llll}
\displaystyle \min_{x_{1}\in\Re^{n_1}}& c_{1}^\top  x_{1} + \check \Q_{2}(x_{1})(\xi_1)\\
\mbox{s.t.}& A_{1}x_{1}=b_{1}\\
& x_{1}\in X_{1}
\end{array}
\right. \nonumber \\ & \equiv \left\{
\begin{array}{llll}
\displaystyle \min_{x_{1},r_{2}\in\Re^{n_1}\times\Re}& c_{1}^\top  x_{1} + r_{2}\\
\mbox{s.t.}&  A_{1}x_{1}=b_{1}\\
 & \beta_{2}^{j\top}(\xi_1) x_{1} + \alpha_2^j(\xi_1) \leq r_{2},\quad j\in J_2\\
 &  x_{1}\in X_{1}\,.
\end{array}
\right. \label{cp1}
\end{align}
The value $\underline{z}$ is a lower bound for the optimal value
of \eqref{mslp}. The computed cutting-plane models $\check \Q_t$,
$t=2,\ldots, T$, and a solution $\bar x_1$ of problem \eqref{cp1}
can be used for constructing an implementable policy as follows.

\paragraph{Forward step.} Given a scenario $\xi=(\xi_1,\ldots,\xi_T) \in
\Xi_1 \times\ldots\times \Xi_T$ (realization of the stochastic
process), decisions $\bar x_t=\bar x_t(\xi_{[t]})$,
$t=1,\ldots,T$, are computed recursively going forward with $\bar
x_1$ being a solution of \eqref{cp1}, and $\bar x_t$ being an
optimal solution of
\begin{equation}
\left\{
\begin{array}{llll}
\displaystyle \min_{x_{t}\in\Re^{n_t}}& c_{t}^\top  x_{t} + \check \Q_{t+1}(x_{t})(\xi_t)\\
\mbox{s.t.}& A_{t}x_{t}=b_{t} - B_{t}\bar x_{t-1}\\
& x_{t}\in X_{t}
\end{array}
\right.  \hspace{-0.1cm} \equiv \left\{
\begin{array}{llll}
\displaystyle \min_{x_{t},r_{t+1}}& c_{t}^\top  x_{t} + r_{t+1}\\
\mbox{s.t.}&  A_{t}x_{t}=b_{t} - B_{t}\bar x_{t-1}\\
 & \beta_{t+1}^{j\top}(\xi_t) x_{t} + \alpha_{t+1}^j(\xi_t) \leq r_{t+1},\quad j \in J_{t+1}\\
 &  x_{t}\in X_{t}\subseteq \Re^{n_t}, r_{t+1} \in \Re \,.
\end{array}
\right. \label{forwt}
\end{equation}
for all $t=2,\ldots,T$, with $\check \Q_{T+1}\equiv 0$. Notice
that $\bar x_t$ is a function of $\bar x_{t-1}$ and
$\xi_t=(c_t,A_t,B_t,b_t)$, i.e., $\bar x(\xi_{[t]})$ is a feasible
and implementable policy for problem~\eqref{mslp} (up to stage
$t$). As a result, the value
\begin{equation}\label{zsup}
\overline{z}=\E\left[\sum_{t=1}^T c_t^\top \bar
x_t(\xi_{[t]})\right]
\end{equation}
is an upper bound for the optimal value of \eqref{mslp} as long as
all the scenarios $\xi \in \Xi$ are considered for computing the
policy. This is the case of nested decomposition. However, the
forward step of the SDDP algorithm consists in taking a sample
$\mathcal{J}$ with $M<N$ scenarios $\xi^j$ of the data process and
computing $\bar x_t(\xi_{[t]}^j)$ and the respective values
$z(\xi^j) = \sum_{t=1}^T c_t^\top \bar x_t(\xi_{[t]}^j)$,
$j=1,\ldots,M$. The sample average $\tilde  z=
\E_{|\mathcal{J}}[z(\xi)]$ and the sample variance $\tilde
\sigma^2=\E_{|\mathcal{J}}[(z(\xi^j)-\tilde z)^2$ are easily
computed. The sample average is an unbiased estimator of the
expectation~\eqref{zsup} (that is an upper bound for the optimal
value of \eqref{mslp}). In the case of subsampling, $\tilde z +
1.96\tilde \sigma/\sqrt{M}$ gives an upper bound for the optimal
value of~\eqref{mslp} with confidence of about $95\%$. As a
result, a possible stopping test for the SDDP algorithm is $\tilde
z + 1.96\tilde \sigma/\sqrt{M} -\underline{z}\leq \epsilon$, for a
given tolerance $\epsilon>0$. We refer to
\cite[Sec.3]{Shapiro_2011} for a discussion on this subject.

\section{On Conditional cuts in SDDP}
\label{sec:condcut}

The assumption that the random data process is stagewise
independent is useful in order to simplify the expressions given
in section \ref{sec:decomp} such as \eqref{cutcoeff}. Essentially
it deletes the additional dependency on the stochastic process of
the cuts. This has a clear advantage in terms of storage, but the
stagewise independence is also unrealistic in many situations. We
will present here an alternative method for dealing with the
conditional expectations that has the following advantages:
\begin{itemize}
\item There is no need to explicitly set up a scenario tree when
initial data is available in the form of a set of scenarios
$\col{(\xi_1^s,...,\xi_T^s)}_{s \in \scr{S}}$.
\item There is no need to increase the dimension of the
state-vector in order to account for stagewise dependency. Ways of handling dependency within SDDP are given in \cite{Guigues_2014,Infanger_Morton_1996}, wherein the specific structure is exploited so as to efficiently ``share" cuts. Similar ideas were successively exploited in \cite{vanAckooij_Malick_2014} but not for SDDP. The methodology laid out in this work does not require making explicit such specific structure.
\end{itemize}
However our underlying assumption is that the random process is
Markovian, i.e., the conditional distribution of $\xi_{t+1}$
knowing $\xi_{[t]}$ is the same as the conditional $\xi_{t+1}$
knowing $\xi_{t}$. We translate this here in the following form:
\begin{equation} \label{eq:xiEvol}
\xi_{t+1} = f_t(\xi_t, \eta_t),
\end{equation}
where the random data or innovation process $\eta_t \in \Re^{p_t}$
is independent of $\xi_{[t]}$ and $f_t$ is a given map defined on
$(\Re^{m_t} \times \Re^{p_t})$ with values in $\Re^{m_{t+1}}$.
This form \eqref{eq:xiEvol} fits general time series models, such
as AR, etc... (see e.g., \cite{Shumway_Stoffer_2005}). Of course exploiting explicitly this structure to compute the conditional expectations, even if the distributions of $\eta_t$ are known is not easy. 

A final approach to represent uncertainty, that we have not explicitly mentioned so far is the use of Markov Chains. We refer to \cite{Mo_Ghelsvik_Grundt_2001,Philpott_deMatos_2012,Valladao_Silva_Poggi_2019} for some references incorporating SDDP and Markov Chains. Finite state Markov chains can of course explicitly be expanded into trees whenever necessary. Therefore in principle similar methodology can be employed. However much like our ``criticism" of scenario trees, it is necessary to ``create" a Markov chain if the initial data is available in the form of a set of scenarios. Since the Markov setting does not fit well the numerical experiments in the paper, we do not provide any comparison with it.

%
%

\subsection{Dealing with conditional expectations - Linear Regression}
\label{subsec:reg}

One of the difficulties in solving problem \eqref{mslp} is closely
tied in with evaluating conditional expectations. Indeed, e.g.,
equation \eqref{cutcoeff} makes it apparent how this introduces an
additional dependency with respect to the situation wherein only
expectations are evaluated. The easiest way of computing
conditional expectations is by setting up a tree or Markov chain
to describe the evolution of $\xi_t$ and employ the transition
probabilities to simply evaluate the conditional expectation.
However, this approach, when only historical data is available,
may lead to modelization errors by introducing an artificial
representation of the underlying process in the form of a tree
when only a set of samples of this process was available. Although
the statistical quality of such a representation can be studied,
we rather suggest to circumvent it altogether by directly
employing a Monte Carlo approach both in the {\it forward} and
{\it backward} steps.

We suppose that in the {\it backward} step the process is sampled
through the Markovian dynamic of equation \eqref{eq:xiEvol} once
and for all with $S$ samples $\xi_t^{k}$ with $k=1,...,S$. This
could be the result of having been able to fit a certain model to
some original data, or having been given a model (e.g., electrical
load could follow the dynamics of the model exhibited in
\cite{Bruhns_Deurveilher_Roy_2005}). Alternatively, these samples
could simply be the historical data supposed to be Markovian. In
the latter case, the data set could just be observed historical
information (e.g., yearly electrical load over the past 30 years
or so). At the end of this preprocessing step, a set of $S$ paths
is available to work with.

It is well known that for any given $t$ the conditional
expectation $\E_{|\xi_{t}}(.)$ is an orthogonal projection in the
space of square integrable  function so we will suppose in the
sequel that all functions $g$ of $\xi = \col{\xi}_{ t >0}$ of
which we want estimate the conditional expectation satisfy $\E[
g(\xi_t)^2] < \infty$ for $t \le T$ : such functions $g$ are said
to be bounded in $L_2$. As a consequence, for any map $g
\fc{\Xi_{t+1}}{\Re}$, the conditional expectation satisfies
$\E_{|\xi_{t}}(g(\xi_{t+1})) = F(\xi_{t})$ for a certain function
$F \fc{\Xi_t}{\Re}$. In order to keep the computational burden
manageable as well as to limit storage requirements, we will
assume that $F$ can be represented as a linear combination of a
given set of base-functions. For instance, one could take
Chebyshev or Legendre polynomials known to be an orthonormal basis
for $L_2$  bounded functions as was originally proposed by
\cite{Longstaff_Schwartz_2001}.

Now for a fixed $t$ and given set $\psi_{t,1},\ldots,
\psi_{t,P_t}$ of base functions the map $F$ can be approximated
with the help of $S$ drawn samples by (linearly) regressing
$\col{g(\xi_{t+1}^s)}_{s=1}^S$ on
$\col{(\psi_{t,1}(\xi_t^s),...,\psi_{t,P_t}(\xi_t^s))}_{s=1}^S$.
This gives the  following approximation:
\begin{equation}
\hat \E^S(g(\xi_{t+1})~|~\xi_{t}) =  \sum_{i=1}^{P_t}
\alpha^*_{t,i}\psi_{t,i}(\xi_t), \label{eq:linreg}
\end{equation}
where $\alpha^*_{t,1},...,\alpha^*_{t,P_t}$ is the optimal
solution to the problem
\begin{equation*}
\min_{\alpha_1,...,\alpha_{P_t}} \sum_{s=1}^{S} \left (
g(\xi_{t+1}^s) - \sum_{i=1}^{P_t} \alpha_{i}\psi_{t,i}(\xi_t^s)
\right )^2.
\end{equation*}

\begin{figure}[!htbp]
\centering
\includegraphics[width=2.5in]{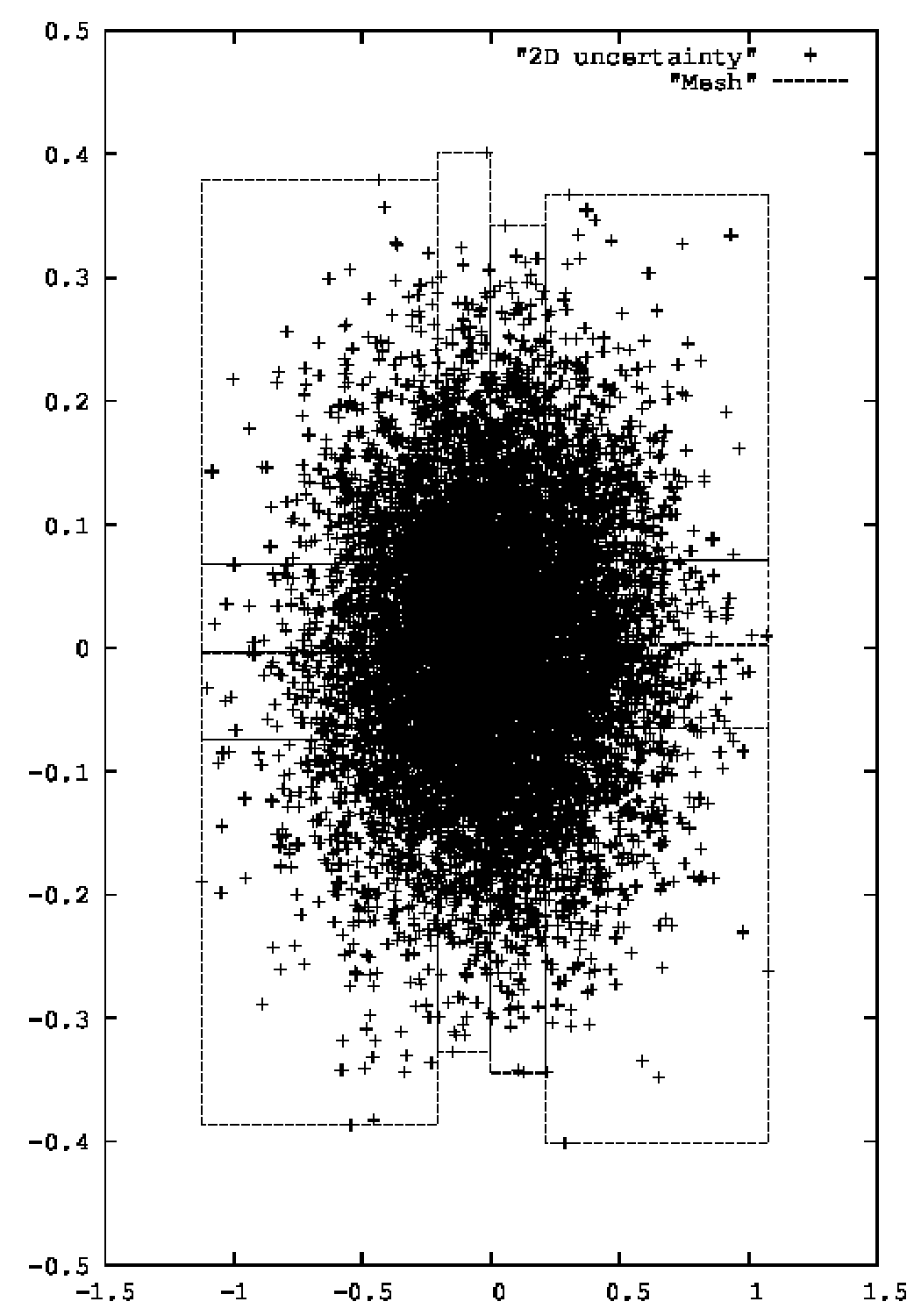}
\caption{ Example of the support of the local basis functions in
dimension 2. Here $4$ intervals are chosen for each direction
giving a total of $L = 16$ hypercubes. \label{fig:mesh}}
\end{figure}

Since the choice of an appropriate well performing basis
(uniformly in $t$) and for a given set of instances of
\eqref{mslp} may be tricky, we will consider instead the adaptive
method developed in \cite{Bouchard_Warin_2012}. It relies on
setting up a class of local base functions. In order to present
the idea, let us denote with $\xi_{t,i}^{s}$ the $i$th coordinate
of sample $s$ of $\xi_t \in \Re^{m_t}$, where $i=1,...,m_t$. Now
taking coordinate-wise extrema over the set of scenarios we define
for each $i=1,...,m_t$, $\bar \xi_{t,i} = \max_{s=1,...,S}
\xi_{t,i}^{s},$  $\underline \xi_{t,i} = \min_{s=1,...,S}
\xi_{t,i}^{s}$. As shown in Figure \ref{fig:mesh}, we now
partition the set $\prod_{i=1}^{m_t} [\underline \xi_{t,i},\bar
\xi_{t,i}]$ in a set of hypercubes $D_{t,\ell}$, $\ell=1,...,L_t$.
For each given $\ell=1,...,L_t$, i.e., on each hypercube, we pick
a family of mappings $\col{\psi_q}_{q \geq 0}$ having support on
$D_{t,\ell}$. Typically the family $\col{\psi_q}_{q \geq 0}$ will
be that of all monomials in $\Re^{m_t}$, implying that $m_t + 1$
coefficients need to be stored for each element $\ell=1,...,L_t$
of the partition. We emphasize that the approximation is
nonconforming in the sense that we do not assure the continuity of
the approximation. However, it has the advantage to be able to fit
any, even discontinuous, function. The total number of degrees of
freedom is equal to $(1+m_t)L_t$ which quality-wise should be
related to the sample size $S$. Hence, in order to avoid
oscillations, the partition is set up so that each element roughly
contains the same number of samples. Using this approximation, the
normal equation is always well conditioned leading to the
possibility to use the Choleski method, which, as shown in
\cite{Bouchard_Warin_2012}, is most efficient for solving the
regression problem. We even note that the problem can be
decomposed in $L_t$ regressions with $m_t \times m_t$ regression
matrices. Using a partial sort algorithm in each direction it is
then possible to get roughly the same number of samples in each
cell as explained in \cite{Bouchard_Warin_2012}. The use of this
partitioning procedure is illustrated on Figure \ref{fig:mesh}.

As explained in the introduction, the numerical results presented
in \cite{Bouchard_Warin_2012} demonstrate the superiority of this
local method with respect to various competing ways of estimating
conditional expectations such as quantization or Malliavin
regression. We note that a constant per mesh approximation can be
used instead of a linear one, giving a far less efficient method
if the functions to approximate are regular or piecewise regular.

\subsection{The SDDP algorithm with conditional cuts and local base functions}
\label{algo:genSDDP}

In this section we provide the description of our variant of the
SDDP algorithm relying on conditional cuts and using local base
functions. A variant using only global base functions is readily
extracted. In the algorithm below, we add the index $s$ to
emphasize the dependency of the data on the scenario. Let us also
note that for a given current trial decision $\bar x_t=\bar
x_t(\xi_{[t]})$ at stage $t$, that this can be associated with a
given mesh $D_{t,h(t)}$, for $h(t) \in \{1,...,L_t \}$. The value
$h(t)$ is such that $\xi_t \in D_{t,h(t)}$. We will emphasize this
by speaking of the trial decision $(\bar x_t, h(t))$.

\paragraph{Backward step.}

Let $(\bar x_t, h(t))$ be a trial decision at stage $t = 1,\ldots,
T-1$ and $\check{\Q}_t$ be a current approximation of the
cost-to-go function $\Q_t$, $t = 2,\ldots, T$, given by the
maximum of a collection of (functional) cutting planes. At stage
$t = T$ the following problem is solved for all $s$ such that
$\xi_{T-1}^{(s)} \in D_{T-1,h(T-1)}$
\begin{equation}\label{cpTEsCond}
\underline{Q}_T(\bar x_{T-1},\xi_T^{(s)})= \left\{
\begin{array}{llll}
\displaystyle\min_{x_T\in\Re^{n_T}} & (c_T^{(s)})^\top  x_T\\
\mbox{s.t.} & A_T^{(s)} x_T=b_T^{(s)} - B_T^{(s)} \bar x_{T-1} \\
& x_T\in X_T
\end{array}
\right.
\end{equation}
Let $\bar \pi_T^{(s)}$ be an optimal dual solution of problem
\eqref{cpTEsCond}. Now let $\alpha_T(y)$, $\beta_T(y)$ be the
numerically approximated conditional expectations of $b_T^\top
\bar \pi_T$ and $B_T^\top \bar \pi_T$ respectively while employing
the estimates of \eqref{eq:linreg}. We will denote this as
$\alpha_T(y):= \hat \E^S[b_T^\top \bar \pi_T | \xi_{T-1}=y]$ and $
\beta_T(y) := -\hat \E^S[B_T^\top \bar \pi_T | \xi_{T-1}=y]$.

Now define the linearization
\begin{align*}
q_T(x_{T-1})(\xi_{T-1}) &:= \beta_T^{\top}(\xi_{T-1}) x_{T-1} +
\alpha_T(\xi_{T-1}) \\ &= \hat \E^S\left[\underline{Q}_T(\bar
x_{T-1},\xi_T)|\xi_{T-1}\right] +
\inner{\beta_T(\xi_{T-1})}{x_{T-1} - \bar x_{T-1}}.
\end{align*}

This linearization is used to update the current cutting plane
model as follows:
\begin{equation*}
\check\Q_T(x_{T-1})(\xi_{T-1}):= \max\{\check \Q_T(x_{T-1})
(\xi_{T-1}),q_T(x_{T-1})(\xi_{T-1})\}.
\end{equation*}
We note that for $\xi \notin D_{T-1, h(T-1)}$,
$q_T(x_{T-1})(\xi)=0$ so that the model is only locally updated.
We will emphasize this by letting $J_{t+1}(D_{t,h(t)})$ be the
index set of linearization added while exploring mesh
$D_{t,h(t)}$.


The updated model $\check \Q_T$ is then used at stage $T-1$, and
the following problem needs to be solved for all $t=T-1,\ldots, 2$
and every $s$ such that $\xi_{t-1}^{(s)} \in D_{t-1,h(t-1)}$:
\begin{align}\label{cpt-1EsCond}
\underline{Q}_{t}(\bar x_{t-1},\xi_{t}^{(s)})=& \left\{
\begin{array}{llll}
\displaystyle \min_{x_{t}\in\Re^{n_t}}& (c_{t}^{(s)})^\top  x_{t} + \check \Q_{t+1}(x_{t})(\xi^{(s)}_t)\\
\mbox{s.t.}& A_{t}^{(s)} x_{t}=b^{(s)}_{t} - B_{t}^{(s)} \bar x_{t-1}\\
& x_{t}\in X_{t}
\end{array}
\right. \\
\quad \equiv \quad &\left\{
\begin{array}{llll}
\displaystyle \min_{x_{t},r_{t+1}\in\Re^{n_t}\times\Re}& (c_{t}^{(s)})^\top  x_{t} + r^{(s)}_{t+1}\\
\mbox{s.t.}&  A_{t}^{(s)} x_{t}=b_{t}^{(s)} - B_{t}^{(s)} \bar x_{t-1}\\
 & \beta_{t+1}^{j\top}(\xi_t^{(s)}) x_{t} + \alpha_{t+1}^j(\xi_t^{(s)}) \leq r_{t+1}^{(s)},\quad j \in J_{t+1}(D_{t,h(t)})\\
 &  x_{t}\in X_{t}\,.
\end{array}
\right.
\end{align}
Let $\bar \pi_t^{(s)}$ and $\big(\bar  \rho_t^{j} \big)^{(s)}$ be
optimal Lagrange multipliers associated with the constraints
$A_{t}^{(s)} x_{t}=b_{t}^{(s)} - B_{t}^{(s)} \bar x_{t-1}$ and
$\beta_{t+1}^{j\top}(\xi_t^{(s)}) x_{t} +
\alpha_{t+1}^j(\xi_t^{(s)}) \leq r_{t+1}^{(s)}$, respectively.
Then the linearization defined for $\xi_{t-1} \in D_{t-1,h(t-1)}$:
\begin{align*}
q_t(x_{t-1})(\xi_{t-1}) &:= \beta_t^{\top}(\xi_{t-1}) x_{t-1} +
\alpha_t(\xi_{t-1}) \\ &= \hat \E^S[\underline{Q}_{t}(\bar
x_{t-1},\xi_{t}) | \xi_{t-1}] + \inner{\beta_t(\xi_{t-1})}{x_{t-1}
- \bar x_{t-1}}
\end{align*}
is constructed with
 \begin{equation}\label{cutcoeffEsCond}
 \alpha_t(\xi_{t-1}):= \hat \E^S[b_t^\top \bar \pi_t + \sum_{j\in J_{t+1}} \alpha_{t+1}^j \bar \rho_t^j | \xi_{t-1}] \quad\mbox{and}\quad
 \beta_t(\xi_{t-1}):= -\hat \E^S[B_t^\top \bar \pi_t |  \xi_{t-1}].
 \end{equation}
Once the above linearization is computed, the cutting-plane model
at stage $t$ is updated

At the first stage, the following problem is solved :
\begin{align}
\underline{z} &=\left\{
\begin{array}{llll}
\displaystyle \min_{x_{1}\in\Re^{n_1}}& (c_{1})^\top  x_{1} + \check \Q_{2}(x_{1})(\xi_1)\\
\mbox{s.t.}& A_{1}x_{1}=b_{1}\\
& x_{1}\in X_{1}
\end{array}
\right. \nonumber \\ &\equiv \left\{
\begin{array}{llll}
\displaystyle \min_{x_{1},r_{2}\in\Re^{n_1}\times\Re}& (c_{1})^\top  x_{1} + r_{2}^{(m)}\\
\mbox{s.t.}&  A_{1}x_{1}=b_{1}\\
 & \beta_{2}^{j\top}(\xi_1) x_{1} + \alpha_2^j(\xi_1) \leq r_{2},\quad j \in J_2\\
 &  x_{1}\in X_{1}\,.
\end{array}
\right. \label{cp1EsCond}
\end{align}
In particular since $\xi_1$ is deterministic $D_{1,1}$ is the
(only) degenerate mesh consisting of $\xi_1$.

\paragraph{Forward step.}
The forward step is essentially unaltered. For a given scenario
$\xi=(\xi_1,\ldots,\xi_T) \in \Xi_1 \times\ldots\times \Xi_T$,
first calculate for each $t$, $h(t) \in \{ 1,...,L_T\}$ such that
$\xi_t \in D_{t,h(t)}$. The decisions $\bar x_t = \bar
x_t(\xi_{[t]})$, $t=1,\ldots,T$, are computed recursively going
forward with $\bar x_1$ being a solution of \eqref{cp1EsCond}, and
$\bar x_t$ being an optimal solution of
\begin{equation}
\left.
\begin{array}{llll}
\displaystyle \min_{x_{t}}\hspace{-0.1cm}& c_{t}^\top  x_{t} + \check \Q_{t+1}(x_{t})(\xi_t)\\
\mbox{s.t.}& A_{t}x_{t}=b_{t} - B_{t}\bar x_{t-1}\\
& x_{t}\in X_{t}\subseteq \Re^{n_t}
\end{array}
\right. \hspace{-0.15cm} \equiv \left\{
\begin{array}{llll}
\displaystyle \min_{x_{t},r_{t+1}}\hspace{-0.15cm}& c_{t}^\top  x_{t} + r_{t+1}\\
\mbox{s.t.}&  A_{t}x_{t}=b_{t} - B_{t}\bar x_{t-1}\\
 & \beta_{t+1}^{j\top}(\xi_t) x_{t} + \alpha_{t+1}^j(\xi_t) \leq r_{t+1},\; j \in J_{t+1}(D_{t,h(t)}),\\
 &  x_{t}\in X_{t}\subseteq \Re^{n_t}, r_{t+1}\in\Re \,.
\end{array}
\right. \label{forwtEsCond}
\end{equation}
for all $t=2,\ldots,T$, with $\check \Q_{T+1}\equiv 0$. The trial
for the next backward recursion is defined as $(\bar x_t, h(t))$
so that in the next backward step cuts at $\bar x_t$ will only be
generated for the visited meshes $D_{t,h(t)}$. We choose to pick
the same stopping criteria as in any usual implementation of SDDP
(see \cite[Remark 1]{Shapiro_2011} for a discussion on this
matter). The advantage of only adding cuts for visited meshes is
that this mitigates the additional burden of storing additional
coefficients. Moreover this cut may be only relevant locally for a
given level of $\xi_t$.

Using the classical SDDP method when the random data process is
stagewise independent, the number of samples used in the forward
path is generally chosen equal to one. However, due to accumulated
errors during the backward resolution the $\check\Q_t$ estimation
based on cutting planes can be far from the corresponding true
supporting planes. In the case of SDDP with conditional cuts, we
want to add one cut for each mesh $D_{t,h(t)}$. Because at each
date $t$, the probability for $\xi_t$ to be in a given mesh is
$\frac{1}{L_t}$, it is numerically better to take a number of
samples equal to $L_t$.

\begin{remark}
For a given accuracy, in the case where a constant per mesh
approximation is used to estimate the conditional expectations,
the number of Monte Carlo trajectories $S$ and the number of
meshes $L_t$ taken are found to be higher than in the linear
approximation. So both the backward resolution and the forward
resolution are more costly with the constant per mesh
approximation than with the linear approximation. We therefore
suggest to use linear per mesh approximation. We refer the reader to Chapter 3 \cite{Friedman_Hastie_Tibshirani_2001} for more information on linear methods for regression and to \cite{Gobet_Lemor_Warin_2005,Lemor_Gobet_Warin_2006} for a cases wherein Monte Carlo error estimates can be given while assuming constant mesh size;
\end{remark}

\section{Numerical experiments}
\label{sec:numerical}

The SDDP method with or without conditional cuts have been
implemented in StOpt \cite{Warin_StOpt} an open source library
containing tools to solve optimization problems both in continuous
or discrete time. In particular, inventory / stock problems (such
as cascaded reservoir management problems, e.g.,
\cite{deMatos_Morton_Finardi_2016,Guigues_Sagastizabal_2012,vanAckooij_Henrion_Moller_Zorgati_2011b}
or gas storage, e.g.,
\cite{Henaff_Laachir_Russo_2013,Carmona_Ludkovski_2010}) can be
optimized using either the Dynamic Programming method or the
Stochastic Dual Dynamic Programming method. The subsequent
experiments were carried out using this library. All computations
are achieved on a cluster composed of Intel Xeon CPU E5-2680 v4
processors whereas all linear programs \eqref{cpt-1EsCond},
\eqref{cp1EsCond} and \eqref{forwtEsCond} are solved using COIN
CLP version 1.15.10.

\subsection{Valuing storage facing the market}
\label{sec:firstTestCase}

We begin with a test case involving gas storage in one dimension
and we will extend it artificially in higher dimensions. All
experiments are achieved with  two processors and 14 cores each.

\subsubsection{The initial test case in one dimension.}

We suppose that we want to value a storage facility where the
owner has the possibility to inject gas from the market or
withdraw it to the market in order to maximize his expected
earnings. The price model is a classical HJM model used for
example in \cite{Warin_2012}:
\begin{equation}
\label{eq:price}
  \frac{dF(T,t)}{F(T,t)} = \sigma e^{-\alpha(T-t)} dW_t
\end{equation}
where $\sigma$ is the volatility of the model, $\alpha$ the mean
reverting, $W$ a Brownian motion on a probability space
$(\Omega,{\cal F},\mathbb{P})$ endowed with the natural (completed
and right-continuous) filtration $\mathbb{F}=\col{{\cal
F}_t}_{t\le T}$ generated by  $W$ up to some fixed time horizon
$T>0$. The spot price is naturally defined as $S_t := F(t,t)$.

We suppose that each date, $t_i= i \Delta T$  with $\Delta T =
\frac{T}{N}$,  $i=0,..., N$, the gas volume $C_i$ has to satisfy
the constraint  $ 0 \le C_i \le C_{\tt max}$. Noting $a_{\tt
in}>0$ the maximal injection during  $\Delta T$,  $a_{\tt out}>0$
the maximal withdrawal from the storage that we suppose to be
constant for simplification purposes, the gain function associated
with a command $x_i \in [-a_{\tt out}, a_{\tt in}]$  at date $ i
\Delta t$ is then
\begin{equation*}
\phi_i(x_i)= -S_{i \Delta T} x_i,
\end{equation*}
where the price $S_{i \Delta T}$ is supposed to be constant
between $t_i$ and $t_{i+1}$. Besides the associated flow equation
is  given by
\begin{equation*}
C_{i+1}= C_i + x.
\end{equation*}
The gain function associated with a strategy $x=\col{x_i
}_{i=0,N}$ is a function depending on the initial storage value
$C_0$, the initial spot price $S_0$:
\begin{equation*}
J(C_0,x,S_0) =  \sum_{i=0}^N \phi_i(x_i),
\end{equation*}
and the optimization problem can be written as
\begin{equation*}
J^*(C_0,S_0) = \sup_{x \in \mathop{U}} J(C_0,x,S_0),
\end{equation*}
where $\mathop{U}$ is the set of non anticipative feasible
strategies.

In this test case, we suppose that we want to optimize the assets
on a one year $(T=1)$ time horizon with annual values $\sigma =
0.5$, $\alpha=0.5$. The decision is taken every week so $N=52$
and we pick the following initial forward curve:
\begin{equation}
\label{forwardPrice} F(0,i \Delta T) = 50 +  \sin(4 \pi
\frac{i}{N}).
\end{equation}
The maximal storage capacity is $C_{\tt max}=200000$ energy units
and the injection and withdrawal rates are $a_{\tt in}= a_{\tt out}=40000$,
 energy units . The ratios
$\frac{C_{\tt max}}{a_{\tt in}}$, $\frac{C_{\tt max}}{a_{\tt
out}}$ indicate a rather slow storage (see \cite{Warin_StOpt} for
different examples of storages). At date $0$ the storage facility
is empty.
This case has a very stochastic effect. The forward curve being very flat,
most of the gain is linked to the volatility of the process. So the resolution is rather difficult.

\begin{figure}[!htbp]
\centering
\subfigure[Price]{
\includegraphics[width=2.5in]{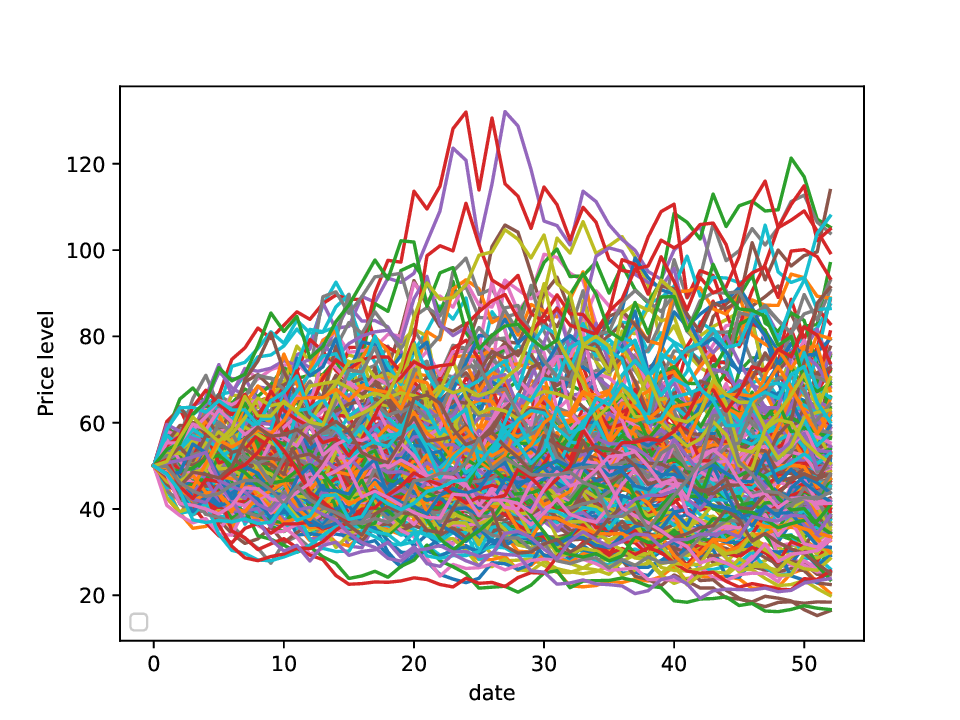}
}
\subfigure[Stock]{
\includegraphics[width=2.5in]{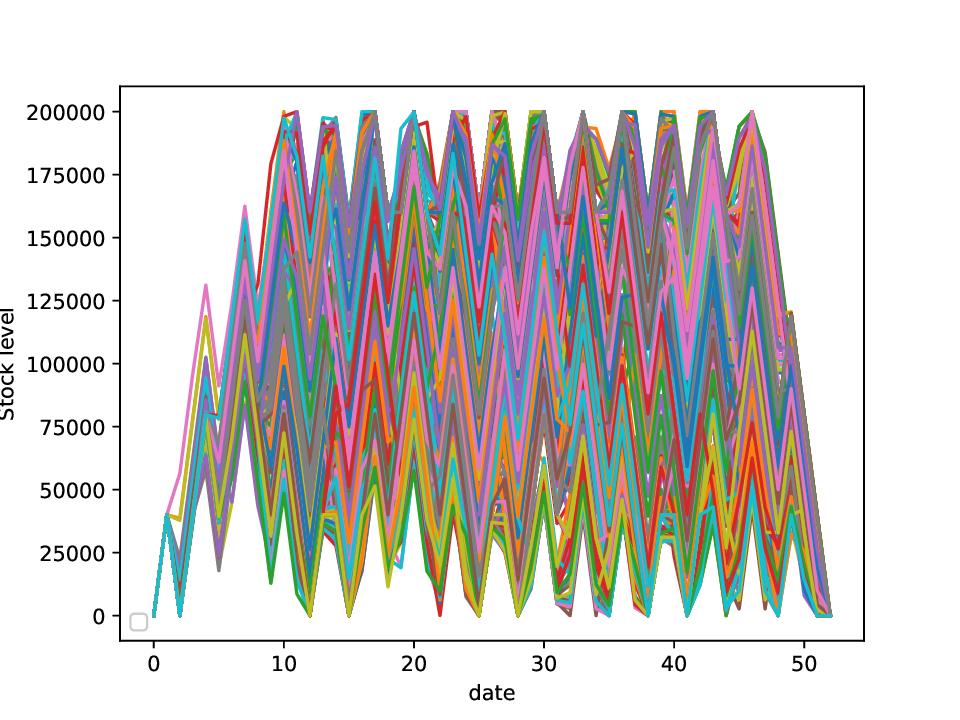}
}
\caption{\label{fig:spotAndTraj} 200 Spot prices trajectories and stock trajectories associated.}
\end{figure}

In table \ref{tab:mc}, we give:
\begin{itemize}
\item the storage valuation when no uncertainty is present (volatility equal to 0),
\item the solution obtained using the dynamic programming method employing the adaptive mesh of \cite{Bouchard_Warin_2012} (DP) and using the assumption that the control is bang-bang (i.e., either equal
to $a_{\tt in}$ or $a_{\tt out}$) \cite{Barrera-Esteve_Bergeret_Dossal_Gobet_Meziou_Munos_Reboul-Salze_2006},
\item the solution obtained using dynamic programming and cuts to estimate the Bellman value (DPC). Regressions using cuts are achieved as exposed in section \ref{subsec:reg} and transition problems are solved using the COIN CLP solver distributing the LPs on the 28
cores.
\end{itemize}
The number of grid points is taken equal to $200$. For DP and DPC results we use $5$ meshes and $20000$ trajectories giving optimization results in row ``Optim Value''. \\
Results obtained using the optimal control in a simulation part, we obtained the results in row `` Simul Value'' using $40000$ trajectories.
\begin{table}
  \centering
  \begin{tabular}{|c|c|c|c|} \hline
    &  Deterministic &  DP & DPC     \\ \hline
    Optim Value&  1.280  &  1.436  &  1.400  \\ \hline
    Simul Value&  1.280  &  1.432  &  1.456      \\ \hline
  \end{tabular}
  \caption{\label{tab:mc} Gas storage valuation in millions of numeraire: deterministic and Monte Carlo solution}
\end{table}
We notice that the deterministic value and the stochastic value are very different, and because prices are martingale, it shows that the problem is  a really stochastic one: optimal trajectories  in simulation and spot prices associated are given on figure \ref{fig:spotAndTraj}. We notice the strong volatility, the mean reversion of the prices  and the diffusion of the stock trajectories.\\
The DCP column permits to show the accuracy of the calculation of the conditional cuts.\\
Another popular way to value storage consist in using scenario trees. We use a trinomial tree \cite{Hull_White_1990} to estimate the solution,with a time step equal to $1/52$. 
In table \ref{tab:valTreeDP}, we give the result obtained by dynamic programming methods using the bang-bang assumption (DP) and solving a LP problem so using cuts associated with tree to estimate the Bellman value (DPC).
\begin{table}
  \centering
  \begin{tabular}{|c|c|c|} \hline
    &   DP & DPC     \\ \hline
    Optim Value&  1.478  &  1.478    \\ \hline
    Simul Value&  1.486 &   1.486       \\ \hline
  \end{tabular}
  \caption{\label{tab:valTreeDP} Gas storage valuation in millions of numeraire by tree methods}
\end{table}
Results obtained by the DP and DPC methods are the same. With the tree method the solution appears to be a little bit higher than with the regression methods. Notice that in the simulation part, uncertainties are sampled on the tree so an error due to the discretization of the tree is present while in the Monte Carlo method with conditional cuts, no discretization error is present. Moreover the approach with a tree is only accurate in very low dimension, as shown in \cite{Bouchard_Warin_2012}. \\

\begin{figure}[!htbp]
\centering
\subfigure[All iterations]{
\includegraphics[width=2.5in]{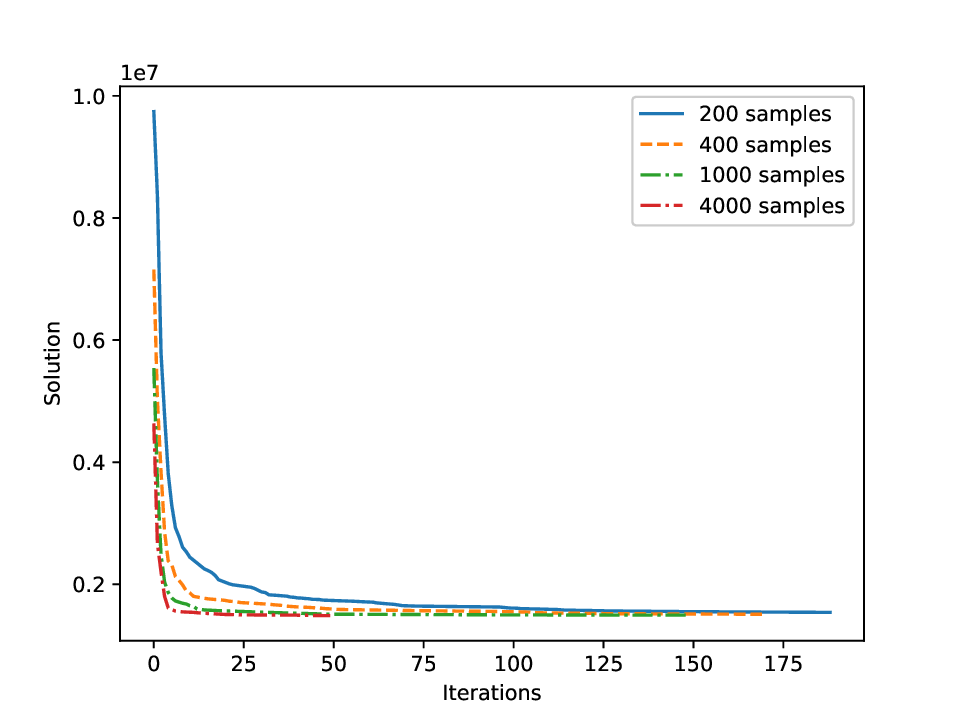}
}
\subfigure[Zoom]{
\includegraphics[width=2.5in]{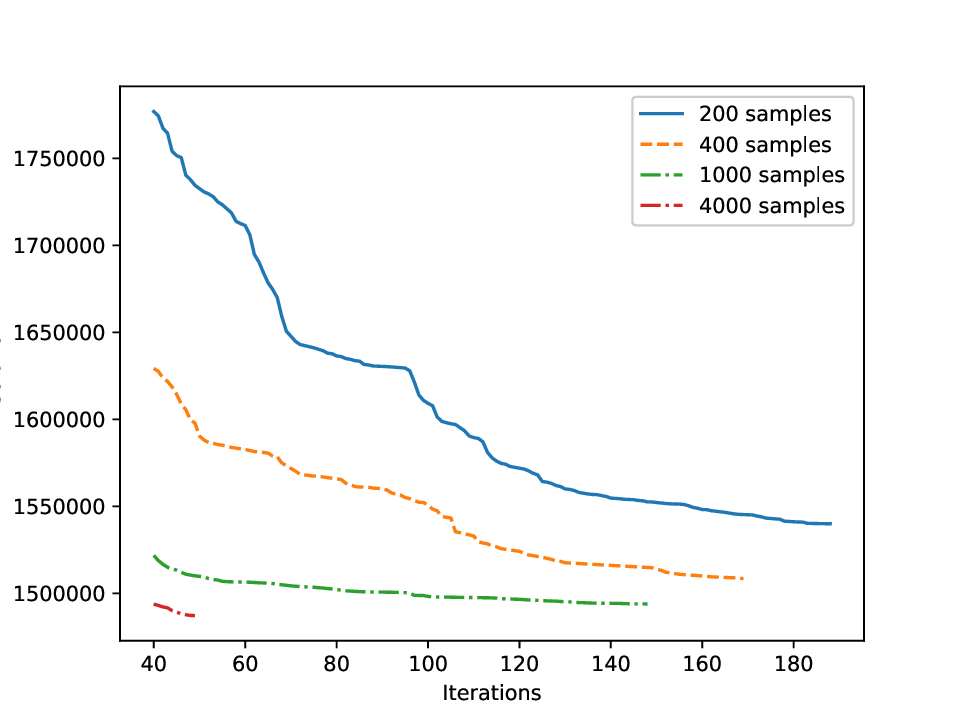}
}
\caption{\label{fig:TreeSDDP} Backward value with trees during iterations depending on the number of samples in forward mode.}
\end{figure}

Next we use the SDDP method with both uncertainties discretized with a trinomial tree and the Monte Carlo method with conditional cuts described in
section \ref{algo:genSDDP}.
The parameters taken are the following:
\begin{itemize}
    \item For the SDDP with a tree, the time step is still $1/52$ and we take a number of samples in the forward part equal to $200$ to $4000$.
    \item For the SDDP with conditional cuts, we take $20000$  trajectories io backward resolution and let the number of mesh vary from $4$ to $6$. The number of samples used in the forward part is taken equal to $28$.
\end{itemize}
While checking the accuracy of the solution we take $40000$ trajectories in a forward part.
For both methods, when a sample gives a new point to explore in the next backward resolution, the cut is only added to the corresponding node in the tree or corresponding cell in the regression method. Moreover no pruning method is used.\\
Results with the tree method are given on figure \ref{fig:TreeSDDP}.
Using many samples, the backward estimation with trees converges far more quickly. Off course on figure \ref{fig:TreeSDDP}, as the number of samples in forward increases we limit the number of iterations as the cost for the resolution explodes. As the number of samples used in forward mode decreases, the convergence rate diminishes dramatically.
It is interesting to note that even if the backward resolution is not very accurate, the forward mode gives good results with not too many iterations. For example, using $1000$ samples in the forward mode, we get a forward value of $1.456e06$ at iteration $50$.\\

\begin{figure}[!htbp]
\centering
\subfigure[All iterations]{
\includegraphics[width=2.5in]{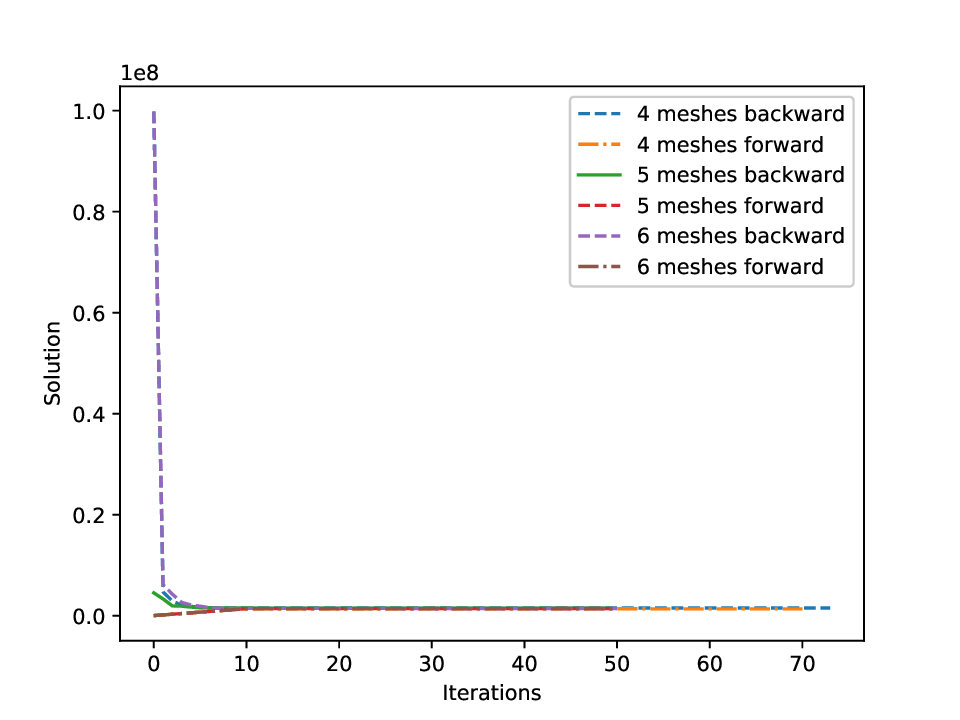}
}
\subfigure[Zoom]{
\includegraphics[width=2.5in]{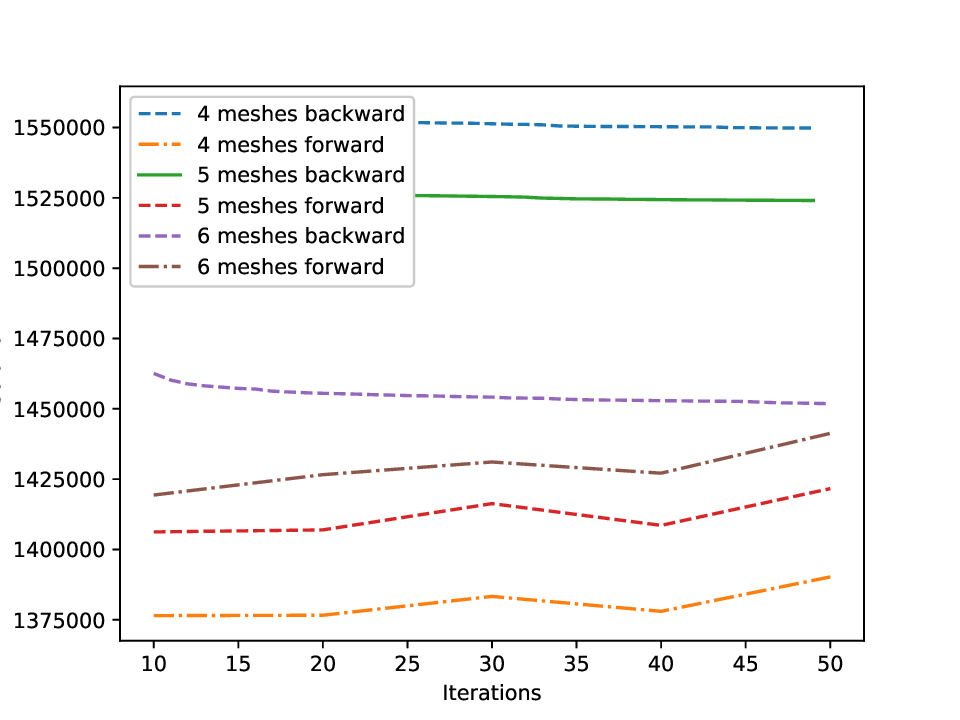}
}
\caption{\label{fig:CondSDDP} Backward value using regressions during iterations depending on the number of meshes.}
\end{figure}

Using regressions, we obtain the results illustrated on figure \ref{fig:CondSDDP}.
Using too few meshes, the conditional cuts are not very well estimated : with 4 meshes  the interval given by the forward values and the backward value is roughly $[1.36e6,1.54e6]$.
Using 5 meshes we get an interval $[1.42e6,1.524e6]$, and using $6$ meshes, we get a tighter estimation $[1.43e6,1.45e6]$.\\
It is to be noticed that when we use to many meshes for a given number of trajectories used in the backward resolution, the backward value can be below the forward one : in this case either we decrease the number of meshes or increase the number of trajectories used in backward  for regressions.\\
Comparing the cost between an SDDP using trees and an SDDP using regression is not very easy. Tighter bounds between forward and backward values are easier to get with regression when the problem is strongly stochastic. In order to have very good values with trees, we have to use many samples in forward mode to explore rarely reached nodes. This drawback is compensated by the fact that each new node explored only leads to three more LP to solve in backward mode (number of branches of the trinomial tree) at each date.

\subsubsection{An artificial test case in higher dimension}

In order to test how the method scales as a function of the
dimension, we take the same data supposing that we dispose of $n$
similar gas storages that we want to value jointly. Each
transition problem is then written as a $n$ dimensional Linear
Program solved with COIN CLP. Of course, in this case, the value
of $n$ storages is simply $n$ times the value of a single storage.

\begin{figure}[!htbp]
\centering
\subfigure[All iterations]{
\includegraphics[width=2.5in]{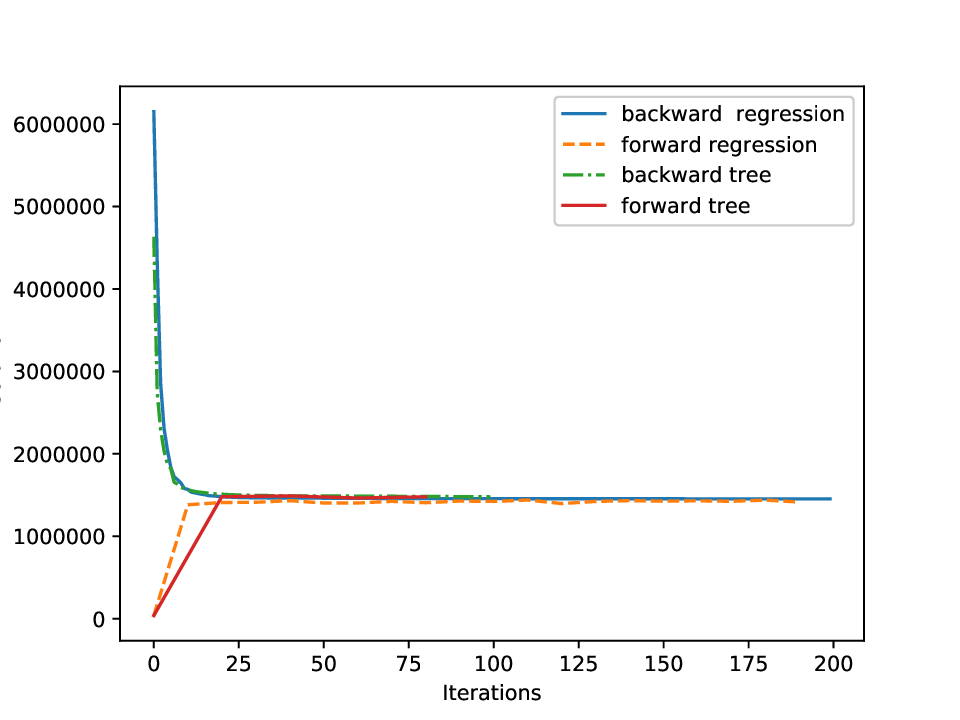}
}
\subfigure[Zoom]{
\includegraphics[width=2.5in]{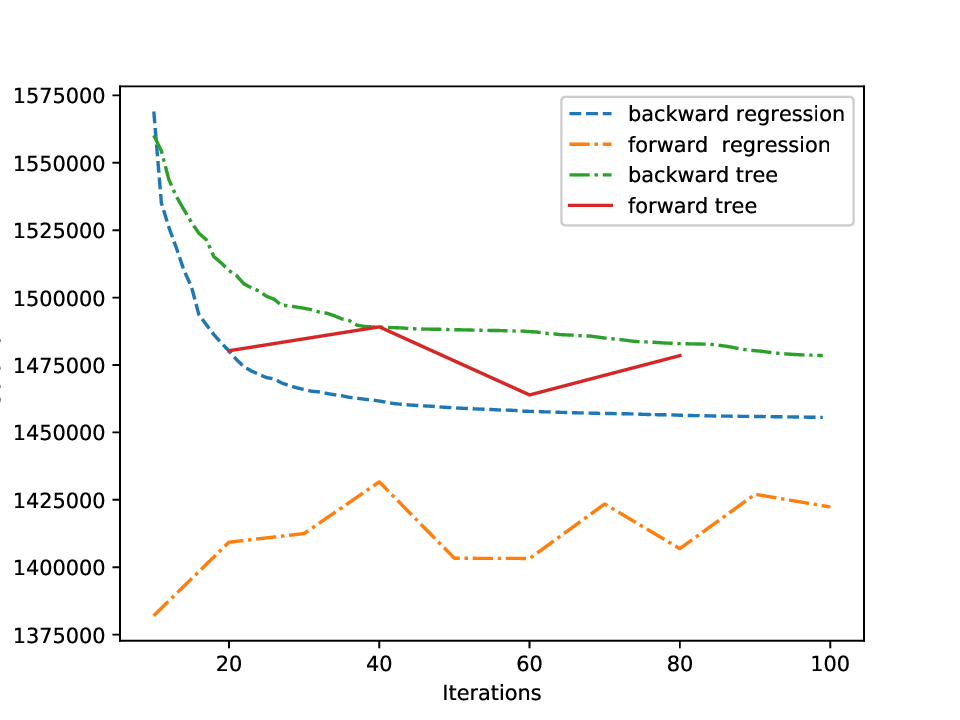}
}
\caption{\label{fig:storage5D} Backward and forward values using regressions (6 meshes) and trees during iterations in 5 dimensions.}
\end{figure}

\begin{figure}[!htbp]
\centering
\subfigure[All iterations]{
\includegraphics[width=2.5in]{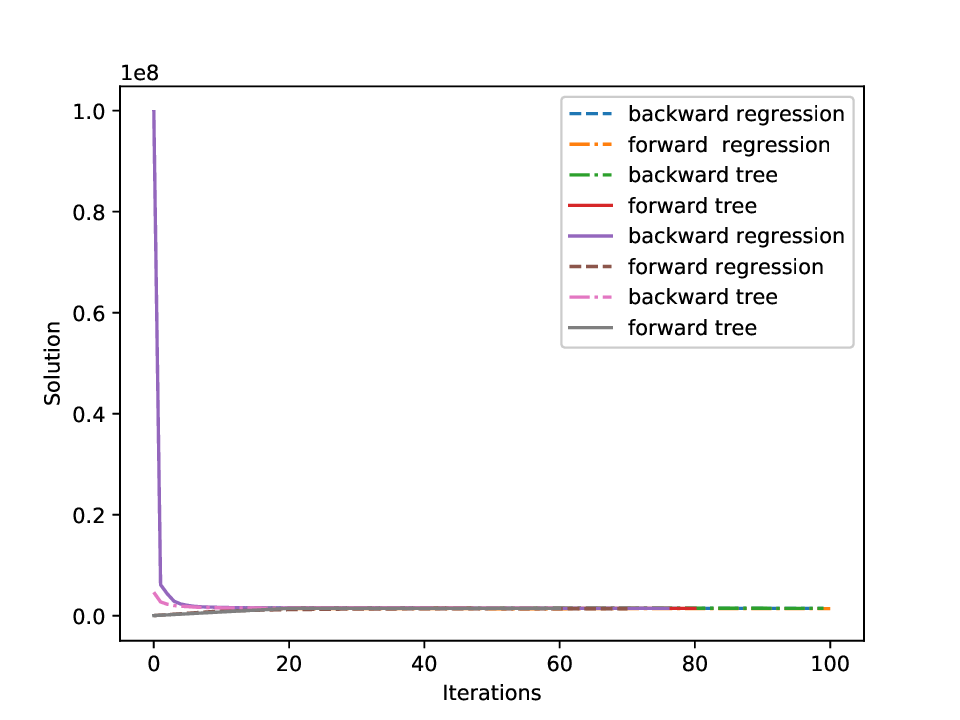}
}
\subfigure[Zoom]{
\includegraphics[width=2.5in]{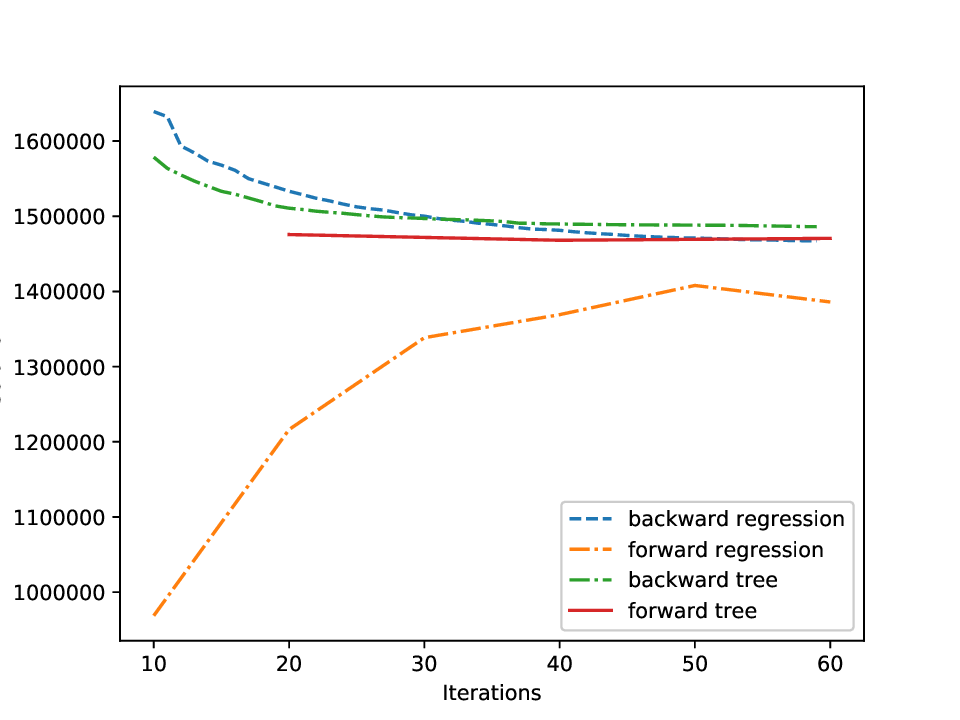}
}
\caption{\label{fig:storage10D} Backward and forward values using regressions (6 meshes) and trees during iterations in 10 dimensions.}
\end{figure}

On figure \ref{fig:storage5D}, \ref{fig:storage10D}, we plot the convergence of
the SDDP scheme for $n=5$ and $n=10$ by giving the value of a
single storage estimated as the value of $n$ storages with the
SDDP algorithm divided by $n$. A forward estimation of the solution is estimated with $40000$ trajectories every 10 SDDP iterations for the regression method and every 20 SDDP iterations for the tree method.  The number of sample used in forward mode is taken equal to $28$ for regressions and $4000$ with trees. As shown in the zoom on the figures, the standard estimation of the solution with Monte Carlo both for tree and regression methods remains rather high with respect to the solution value: oscillations in forward mode are present although the number of trajectories is quite high.

We note that a direct application of dynamic programming
in this situation is simply no longer possible, whereas (our
variant of) SDDP scales up very well.\\
As in the one dimensional case, the value obtained by trees is a little bit higher with the tree method.
In dimension 10, the value obtained with regression is slightly below its one dimensional reference ( a value around $1.40e6$ compared to a value of $1.43e6$ obtained in dimension 1).

\subsection{A simple storage facility problem facing a consumption constraint}
\label{subsec:secondcase}

In this section we compare our conditional cut approach with the
classical SDDP method (i.e., increasing the state-space) for
dealing with uncertainty following an AR(1) model. All experiments
are once again achieved with two processors with 14 cores each.

\subsubsection{A case in one dimension}

\label{subsubsec:load1D} In this test case we suppose that the
price gas is deterministic so is given by the future price
\eqref{forwardPrice}. The characteristics of the storage facility
are the same as in section \ref{sec:firstTestCase} but we suppose
that we have an injection cost per unit of injected gas equal to
$0.1$.

The owner of the storage facility is allowed to buy gas on the
market and to inject it into the storage (paying the injection
cost) but is not allowed to sell it on the market. The gas
withdrawn can only be used to satisfy the consumption $D$
following an AR(1) (e.g., \cite{Shumway_Stoffer_2005}) process
between two dates $t_i= i\frac{T}{N}$ and $t_{i+1} = (i+1)
\frac{T}{N}$:
\begin{equation*}
D_{i+1} - \tilde D_{i+1}= \kappa_d (D_{i}- \tilde D_{i}) +\sigma_d
\epsilon_i,
\end{equation*}
where $\epsilon_i$ is independent white noise with zero mean and
variance 1, $\kappa_d$ is set equal to $0.9$, $\sigma_d=1000$ and
$\tilde D$ is an average consumption rate satisfying:
\begin{equation*}
\tilde D_{i} = 22000 +  7000 \sin(4 \pi \frac{i}{N}).
\end{equation*}

We denote with $x^b_i$ the quantity of gas bought at date $t_i$,
$x^{\tt in}_i \in [0,a_{\tt in}]$ the quantity of gas injected and
$x^{\tt out}_i \in [-a_{\tt out},0]$ the quantity of gas
withdrawn. We thus have the following constraint to satisfy at
each date $t_i$:
\begin{equation*}
x^b_i = D_i + x^{in}_i + x^{out}_i.
\end{equation*}
The following objective function characterizes the cost of a
strategy $x^b=\col{x^b_i}_{i=0,N}$ by the value :
\begin{equation*}
J(C_0,x^b,D_0) = \sum_{i=0}^N -\phi_i(x^b_b),
\end{equation*}
where $D_0$ is the initial consumption value that we suppose to be
equal to $\tilde D_0$ and the optimization problem can be written
as
\begin{equation*}
J^*(C_0,D_0) = \inf_{x \in \mathop{U}} J(C_0,x,D_0).
\end{equation*}

Due to the AR(1) structure of the demand, it is possible to add
$D_i$ to the state vector $x_i$ in the algorithm
\ref{algo:genSDDP} as proposed by \cite{Pereira_Pinto_1991}. We
can then test two approaches:
\begin{itemize}[itemindent=1.8cm, leftmargin=1cm]
\item[-Approach A:] SDDP with conditional cuts as suggested in
section \ref{algo:genSDDP}. Here  $x_i$ is only the stock level
and $\xi_i = ( D_{i})$ is a  dimensional Markov process. All cuts
coefficients are one dimensional functions and have to be
estimated by one dimensional regressions. In the test case we use
$N_c^A=2000$ trajectories in optimization and conditional cuts
coefficients are estimated with an $I_1= 8$ grid. At each step of
the SDDP algorithm we use $8$ forward simulations to discover the
new states used in the next backward step.
\item[-Approach B:] SDDP without  conditional cuts
 but with an augmented state vector. Here $x_i =(x^b_i,
D_i)$, and at each date $t_i$ all cuts coefficients are constant.
This method is related to nested Monte Carlo: a set of $N_i$
simulations is chosen to simulate $\epsilon_t$. So at a given date
in the backward pass, we are left with $N_i$
 LPs to solve for each visited state in the
previous forward simulation.
 It turns out that in order to be able to be able to reach the $0.1\%$ criterium
for difference between the forward and the backward iteration, we
have to take $N_i=2000$. At each step of the SDDP algorithm we use
$1$ forward simulation to discover the new states used in the
next backward step.
\end{itemize}

\begin{figure}[!htbp]
  \centering
  \includegraphics[width=3in]{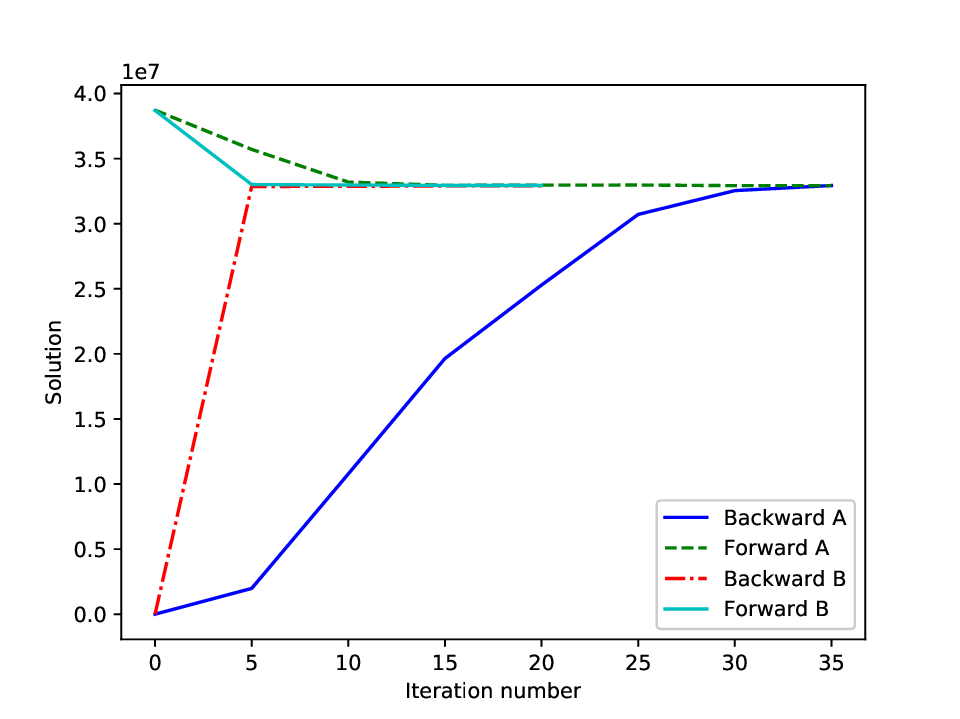}
  \caption{Comparison of optimal values in the backward and
  forward pass for two variants of SDDP for a  load constraint. \label{fig:load1V2}}
\end{figure}

On Figure \ref{fig:load1V2}, we give the results obtained by the
two approaches until convergence of the backward and forward
values to within a given (relative) accuracy of $0.1\%$.
Approach A takes 109 seconds to converge while Approach B takes 66
seconds to converge. On this simple case, Approach B is superior
to Approach A, especially because of a better convergence rate in
the backward step.

\subsubsection{Artificial case in dimension 10}

In order to test how these methods scale up, we pick the test case
of section \ref{subsubsec:load1D} but this time with 10 similar
storage facilities. We also multiply average load as well as
$\sigma_d$ by 10. Consequently the expected value of the optimal
cost is 10 times the value of the case with one storage facility.

\begin{figure}[!htbp]
  \centering
  \includegraphics[width=3in]{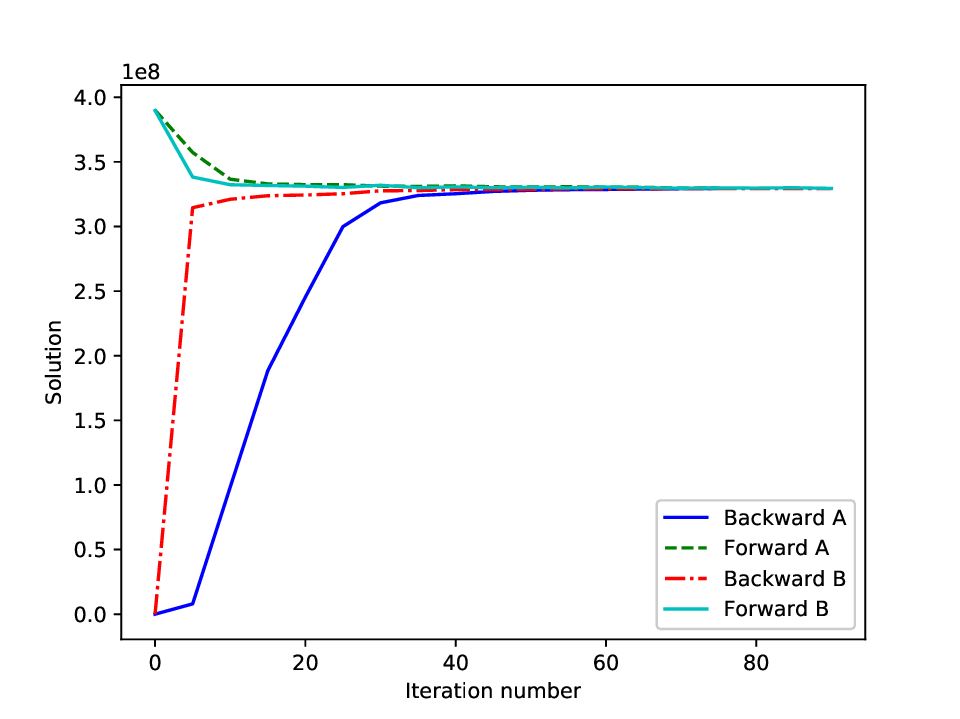}
  \caption{Comparison of optimal values in the backward and
  forward pass for two variants of SDDP for ten gas storage
  facilities with a coupling load constraint. \label{fig:load10V2}}
\end{figure}

On Figure \ref{fig:load10V2}, we give the convergence results
obtained by the two approaches.
In order to be able to get a relative given  accuracy of $0.1\%$
it is necessary to increase $N_i$ to $4000$. Both approaches allow
us to obtain the same value, showing consistency of the suggested
method. Although approach B, obtains a better estimation of the
optimal value earlier on in the algorithm, approach A turns out to
be faster. Indeed, for reaching a 0.1\% asserted gap, Approach A
takes $660$ seconds to converge while Approach B takes $775$
(+17\%) seconds.


\subsection{Valuing a more complex  storage facility when facing a consumption constraint}

In this section we take the test case \ref{subsec:secondcase} and
allow the prices to vary following the dynamic given in section
\ref{sec:firstTestCase}. As in the case described in section
\ref{sec:firstTestCase}, the classical SDDP method cannot deal
with this kind of uncertainty. Note that this is so since price
uncertainty affects the cost coefficients of the transition
problems (e.g., \eqref{cpt-1EsCond}), so that increasing the state
space would make these problems bi-linear (i.e., no-longer
convex). In the two following subsections we show that our method
can deal efficiently with multidimensional uncertainty. All
calculation are achieved on a cluster of 8 processors with 14
cores each.

\subsubsection{A case in one dimension}
\label{subsubsec:gasload1D} Similarly to the sections
\ref{sec:firstTestCase} and \ref{subsec:secondcase} the following
objective function characterizes the cost of a strategy
$x^b=\col{x^b_i}_{i=0,N}$:
\begin{equation*}
J(C_0,x^b,S_0,D_0) = \sum_{i=0}^N -\phi_i(x^b_b).
\end{equation*}
The optimization problem can be written as
\begin{equation*}
J^*(C_0,S_0,D_0) = \inf_{x \in \mathop{U}} J(C_0,x,S_0,D_0).
\end{equation*}

We use the SDDP method with conditional cuts as suggested in
section \ref{algo:genSDDP}. Here  $x_i$ is only the stock level
and $\xi_i = (S_{i \Delta T}, D_{i})$ is a two dimensional Markov
process. All cut coefficients are two dimensional functions and
have to be estimated by two dimensional regressions. In the test
case we use $N_c^A=32000$ trajectories in optimization and
conditional cuts coefficients are estimated with an $(I_1,I_2) =
(8,4)$ grid (see Figure \ref{fig:mesh}). At each step of the SDDP
algorithm we use $32$ forward simulations to discover the new
states used in the next backward step.

\begin{figure}[!htbp]
  \centering
  \includegraphics[width=3in]{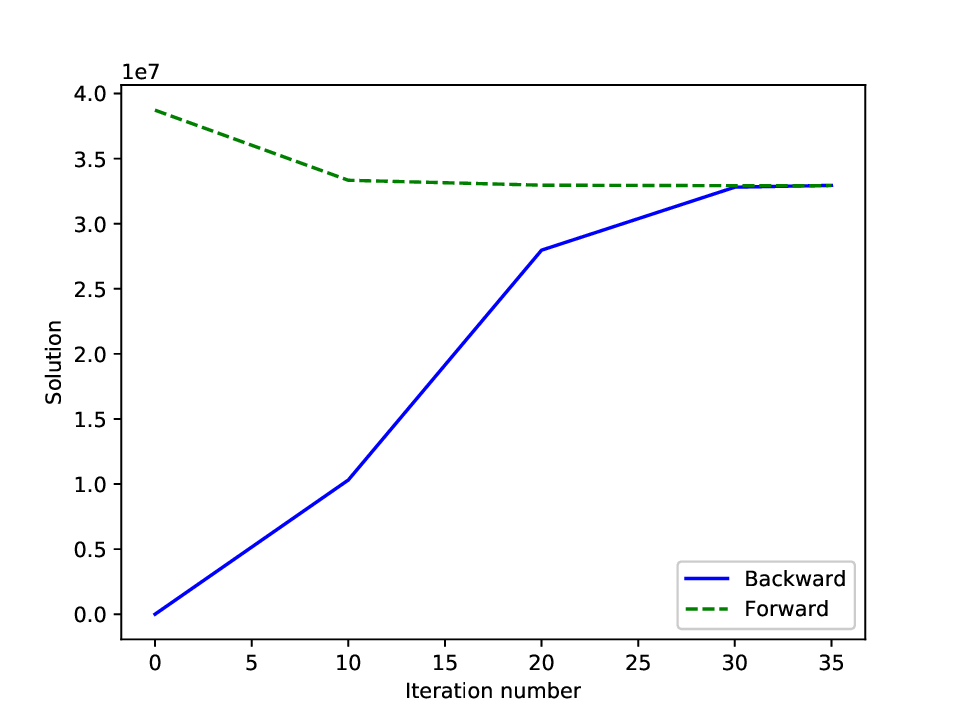}
  \caption{SDDP backward and forward iterations for a  gas storage
  facility with a coupling load constraint. \label{fig:storage1V2}}
\end{figure}

On Figure \ref{fig:storage1V2}, we give the results  until
convergence of the backward and forward values with a given
(relative) accuracy of $0.1\%$. The convergence is achieved in 845
seconds.

\subsubsection{An artificial test case in dimension 10.}

In order to test how the method scales up, we pick the test case
of section \ref{subsubsec:gasload1D} but this time with 10 similar
storage facilities. We also multiply average load as well as
$\sigma_d$ by 10. Consequently the expected value of the optimal
cost is 10 times the value of the case with one storage facility.

\begin{figure}[!htbp]
  \centering
  \includegraphics[width=3in]{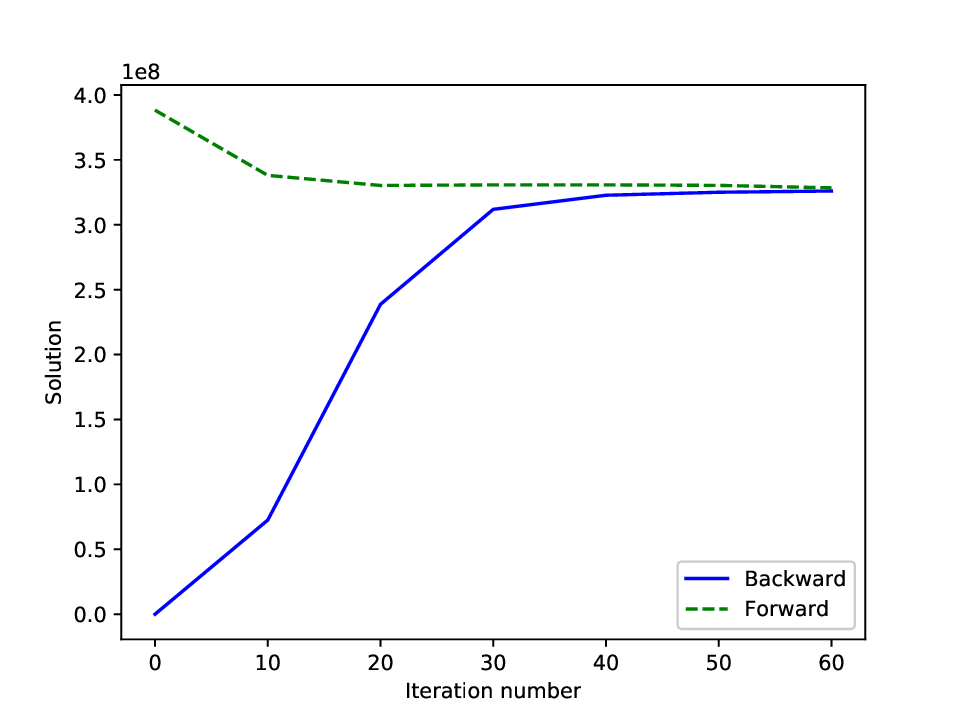}
  \caption{SDDP backward and forward iterations for ten  gas storage
  facilities with a coupling load constraint. \label{fig:storage10V2}}
\end{figure}

On Figure \ref{fig:storage10V2}, we give the convergence results
obtained by the conditional cuts method.
Once again our method converges very efficiently even for a high
dimension problem with convergence reached in 2400 seconds.

\section{Conclusion}

Using conditional cuts we have developed an effective methodology
to value easily high dimensional storages problems. Getting rid of
building trees for the backward part, the method is easy to use
even on a set of initial scenarios, perhaps resulting from
historical data, or any other set of scenarios assumed to be
Markovian. The method can also be employed whatever the elements
of the transition problems impacted by uncertainty. This is not
the case for the classic approach, consisting of increasing the
state space dimension, the use of which is restricted to right
hand side uncertainty.

Moreover, when the dimension of the uncertainties is low and in
the special case of auto regressive models, we have shown that the
use of conditional cuts is competitive with the classical approach
developed in \cite{Pereira_Pinto_1991} consisting of augmenting
the state space to account for past dependency.

At last, because our method is a pure Monte Carlo one, risk
optimization such as CVaR or VaR optimization should be much more
accurate than classical approaches.


\bibliographystyle{plain}
\bibliography{./biblio_vanAckooij_Warin}


\end{document}